\documentclass[12pt]{amsart}

\usepackage{amssymb}
\usepackage{amsthm}
\usepackage{amscd}

\usepackage{pifont}

\usepackage{textcomp}
\usepackage{latexsym}

\usepackage{amsmath}
\usepackage{epic,eepic}
\usepackage[dvipdfm]{graphicx}

\usepackage[all]{xy}

\theoremstyle{plain}
\newtheorem{lemma}{Lemma}[subsection]
\newtheorem{prop}[lemma]{Proposition}
\newtheorem{thm}[lemma]{Theorem}
\newtheorem{cor}[lemma]{Corollary}
\newtheorem{aplemma}{Lemma~A.\hspace{-1.5mm}}
\newtheorem{approp}{Proposition~A.\hspace{-1.5mm}}
\newtheorem{apthm}{Theorem~A.\hspace{-1.5mm}}
\newtheorem{apcor}{Corollary~A.\hspace{-1.5mm}}
\newtheorem{intthm}{Theorem}

\newtheorem{conj}[lemma]{Conjecture}

\theoremstyle{definition}

\newtheorem{rema}[lemma]{Remark}

\newtheorem{remb}{Remark}

\newtheorem{defi}[lemma]{Definition}
\newtheorem{exa}[lemma]{Example}
\newtheorem{aprem}{Remark~A.\hspace{-1.5mm}}
\newtheorem{apdefi}{Definition~A.\hspace{-1.5mm}}
\newcommand{\bde}{\begin{defi}}
\newcommand{\ede}{\end{defi}\vspace{1mm}}
\newcommand{\ble}{\begin{lemma}}
\newcommand{\ele}{\end{lemma}}
\newcommand{\bpr}{\begin{prop}}
\newcommand{\epr}{\end{prop}}
\newcommand{\bt}{\begin{thm}}
\newcommand{\et}{\end{thm}}
\newcommand{\bco}{\begin{cor}}
\newcommand{\eco}{\end{cor}}
\newcommand{\bre}{\begin{rema}}
\newcommand{\ere}{\end{rema}}
\newcommand{\brea}{\begin{rema}}
\newcommand{\erea}{\end{rema}\vspace{1mm}}
\newcommand{\breb}{\begin{remb}}
\newcommand{\ereb}{\end{remb}\vspace{1mm}}
\newcommand{\bex}{\begin{exa}}
\newcommand{\eex}{\end{exa}}
\newcommand{\bpf}{\begin{proof}}
\newcommand{\epf}{\end{proof}\vspace{1mm}}

\newcommand{\bade}{\begin{apdefi}}
\newcommand{\eade}{\end{apdefi}}
\newcommand{\bale}{\begin{aplemma}}
\newcommand{\eale}{\end{aplemma}}
\newcommand{\bapr}{\begin{approp}}
\newcommand{\eapr}{\end{approp}}
\newcommand{\bat}{\begin{apthm}}
\newcommand{\eat}{\end{apthm}}
\newcommand{\baco}{\begin{apcor}}
\newcommand{\eaco}{\end{apcor}}
\newcommand{\bare}{\begin{aprem}}
\newcommand{\eare}{\end{aprem}}

\newcommand{\be}{\begin{enumerate}}
\newcommand{\ee}{\end{enumerate}}
\newcommand{\bcd}{\[\begin{CD}}
\newcommand{\ecd}{\end{CD}\]}
\newcommand{\bit}{\begin{itemize}}
\newcommand{\eit}{\end{itemize}}
\newcommand{\bq}{\begin{quote}}
\newcommand{\eq}{\end{quote}}
\newcommand{\ba}{\begin{array}}
\newcommand{\ea}{\end{array}}

\newcommand{\mcD}{\mathcal{D}}
\newcommand{\mcE}{\mathcal{E}}
\newcommand{\mcF}{\mathcal{F}}
\newcommand{\mcG}{\mathcal{G}}
\newcommand{\mcH}{\mathcal{H}}
\newcommand{\mcI}{\mathcal{I}}

\newcommand{\mcK}{\mathcal{K}}
\newcommand{\mcL}{\mathcal{L}}

\newcommand{\mcN}{\mathcal{N}}
\newcommand{\mcO}{\mathcal{O}}
\newcommand{\mcP}{\mathcal{P}}
\newcommand{\mcQ}{\mathcal{Q}}

\newcommand{\mcS}{\mathcal{S}}
\newcommand{\mcT}{\mathcal{T}}

\newcommand{\mcV}{\mathcal{V}}

\newcommand{\mcZ}{\mathcal{Z}}
\newcommand{\mbA}{\mathbb{A}}

\newcommand{\mbC}{\mathbb{C}}
\newcommand{\mbD}{\mathbb{D}}

\newcommand{\mbG}{\mathbb{G}}
\newcommand{\mbH}{\mathbb{H}}

\newcommand{\mbR}{\mathbb{R}}
\newcommand{\mbS}{\mathbb{S}}

\newcommand{\mbZ}{\mathbb{Z}}

\newcommand{\mfE}{\mathfrak{E}}

\newcommand{\mfM}{\mathfrak{M}}

\newcommand{\mfO}{\mathfrak{O}}

\newcommand{\mfR}{\mathfrak{R}}
\newcommand{\mfS}{\mathfrak{S}}
\newcommand{\mfT}{\mathfrak{T}}
\newcommand{\mfU}{\mathfrak{U}}

\newcommand{\mfc}{\mathfrak{c}}

\newcommand{\mfe}{\mathfrak{e}}

\newcommand{\mfg}{\mathfrak{g}}
\newcommand{\mfh}{\mathfrak{h}}

\newcommand{\mfl}{\mathfrak{l}}

\newcommand{\mfs}{\mathfrak{s}}
\newcommand{\mft}{\mathfrak{t}}

\newcommand{\migi}{\rightarrow}
\newcommand{\longmigi}{\longrightarrow}

\newcommand{\isom}{\stackrel{\sim}{\migi}}

\newcommand{\migiincl}{\hookrightarrow}

\newcommand{\migisurj}{\twoheadrightarrow}

\newcommand{\Dp}{\mfM^{^\text{Zzz...}}}

\newcommand{\G}{G}

\newcommand{\ind}{\circledast}

\newcommand{\mr}{\mathrm}
\newcommand{\hidden}[1]{\,}

\begin{document}

\title[The symplectic nature of dormant indigenous bundles]{The symplectic nature  of  the space \\ of  dormant indigenous bundles \\ on  algebraic curves}
\author{Yasuhiro Wakabayashi}
\date{}
\maketitle
\footnotetext{Y. Wakabayashi: Department of Mathematics, Tokyo Institute of Technology, 2-12-1 Ookayama, Meguro-ku, Tokyo 152-8551, JAPAN;}
\footnotetext{e-mail: {\tt wkbysh@math.titech.ac.jp};}

\begin{abstract}
We study the symplectic nature of the moduli stack classifying  dormant curves over a field $K$ of positive characteristic, i.e.,
proper hyperbolic curves over $K$  equipped with
a dormant indigenous bundle.
The central objects of the present paper are  the following two  Deligne-Mumford stacks.
One is the cotangent bundle 
${^\circledcirc T^{\vee ^\mathrm{Zzz...}}_{g,K}}$
of the moduli stack  ${^\circledcirc \mathfrak{M}^{^\mathrm{Zzz...}}_{g,K}}$ classifying 
ordinary dormant curves over $K$ of genus $g$.
The other is the moduli stack ${^\circledcirc \mathfrak{S}^{^\mathrm{Zzz...}}_{g,K}}$
 classifying ordinary dormant curves over $K$  equipped with an indigenous bundle.
These Deligne-Mumford stacks  admit canonical symplectic structures respectively.
 The  main result  of the present paper  asserts
 that a canonical isomorphism ${^\circledcirc T^{\vee ^\mathrm{Zzz...}}_{g,K}}  \rightarrow  {^\circledcirc \mathfrak{S}^{^\mathrm{Zzz...}}_{g,K}}$ preserves the symplectic structure. 
This result may be thought of as 
 a positive characteristic analogue of the works of S. Kawai  (in the paper entitled ``{\it The symplectic nature of the space of projective connections on Riemann surfaces}"), P. Ar\'{e}s-Gastesi,  I. Biswas, and B. Loustau.
Finally, as its application, we construct  
a Frobenius-constant quantization on the moduli stack 
${^\circledcirc \mathfrak{S}^{^\mathrm{Zzz...}}_{g,K}}$.

\end{abstract}
\tableofcontents 
\section*{Introduction}

\subsection{} \label{s01}
The purpose of  the present paper 
is to  study symplectic geometry of  indigenous bundles 
 in positive characteristic.
Here, recall that the notion of an indigenous bundle was  originally  introduced and studied by R. C. Gunning in the context of 
 Riemann surfaces (cf. ~\cite{G2}).
Roughly speaking,  an indigenous bundle on a connected  compact hyperbolic Riemann surface $X$ is a  projective line bundle (or equivalently, a $\mr{PGL}_2$-torsor) on $X$ together with a  connection and a global section satisfying a certain transversality condition.
It may be thought of as an  {\it algebraic} object  that encodes the ({\it analytic}, i.e., {\it non-algebraic}) uniformization data for the $X$.
Various equivalent mathematical objects,
including  certain kinds of differential operators  between line bundles (related to Schwarzian derivatives),
have been investigated  by many mathematicians.
For example,
   indigenous bundles on $X$ correspond bijectively to  certain  {\it projective structures}
on the underlying topological space $\Sigma$ of $X$.
A projective structure on $\Sigma$ is, by definition, the  equivalence class of an  atlas covered by coordinate charts on $\Sigma$ such that the transition functions are expressed as M\"{o}bius transformations.
(In particular, each projective structure  determines  a unique  Riemann surface structure on $\Sigma$.)
If we are given a projective structure, then the collection of  its transition functions specifies   a $\mr{PGL}_2$-torsor equipped with a connection, which  forms  an indigenous bundle.
In this way, we obtain  a bijective correspondence between the set of projective structures on $\Sigma$ defining the Riemann surface $X$ and the set of (isomorphism classes of) indigenous bundles on $X$.

\vspace{3mm}
\subsection{} \label{s02}

Now, we shall recall the works of S. Kawai, P. Ar\'{e}s-Gastesi, I. Biswas, and B. Loustau (cf. ~\cite{Ka}; ~\cite{AGBi}; ~\cite{AGBi2}; ~\cite{Los}) 
 related to 
  natural symplectic structures on certain moduli spaces of 
projective structures, or equivalently, indigenous bundles.

Let $\Sigma$ be, as above,  a connected orientable closed surface of genus $g>1$.
Write $\mfT^\Sigma$ for the Teichm\"{u}ller space associated with  $\Sigma$, i.e.,
the quotient space
\begin{align}
 \mfT^\Sigma := \mr{Conf}(\Sigma)/\mr{Diff}^0(\Sigma), \notag
 \end{align}
where $\mr{Conf}(\Sigma)$ denotes the space of all conformal structures on $\Sigma$ compatible with the orientation of $\Sigma$, and $\mr{Diff}^0(\Sigma)$ denotes  the group of all diffeomorphisms of $\Sigma$ homotopic to the identity map of $\Sigma$. 
In particular, it is a covering space of the moduli stack  $\mfM_{g,\mbC}^\mr{an}$ classifying connected compact Riemann surfaces of genus $g$.
Also, write
\begin{align} \label{e02}
  \mfS_{\mfT^\Sigma} := \mr{Proj}(\Sigma)/\mr{Diff}^0(\Sigma), 
\end{align}
where $\mr{Proj}(\Sigma)$ denotes the space of all projective structures on $\Sigma$.
It is known that 
the cotangent bundle $T^\vee_{\mfT^\Sigma}$ of $\mfT^\Sigma$ (resp., the quotient space $\mfS_{\mfT^\Sigma}$) admits  a structure of  complex manifold of dimension  $6g-6$, equipped with a  holomorphic symplectic structure $\omega^\mr{can}_{\mfT^\Sigma}$ (resp., $\omega^{\mr{PGL}}_{\mfT^\Sigma}$ (cf. ~\cite{Go})).
Consider an analytic
 section
\begin{align} 
\sigma_{\mr{unif}}  :  \mfT^\Sigma \migi  \mfS_{\mfT^\Sigma}  \notag
 \end{align}
 of the natural projection $\mfS_{\mfT^\Sigma} \migi \mfT^\Sigma$  
arising from the uniformization in the sense of either Bers, Schottky, or Earle (cf. ~\cite{Be}; ~\cite{Bi2}; ~\cite{Ea}). (Note that the Bers uniformization are determined after choosing  a specific point of $\mfT^\Sigma$.)
Because of a natural affine structure on $\mfS_{\mfT^\Sigma}$, 
the section $\sigma_{\mr{unif}}$  may be  extended to   an isomorphism  
\begin{align} 
  \Psi_{\mr{unif}} : T^\vee_{\mfT^\Sigma} \isom \mfS_{\mfT^\Sigma} \notag
  \end{align}
whose restriction  to the zero section $\mfT^\Sigma \migi  T^\vee_{\mfT^\Sigma}$ coincides with $\sigma_{\mr{unif}}$.
It follows from ~\cite{Ka},  Theorem, ~\cite{AGBi}, Theorem 1.1,  ~\cite{AGBi}, Remark 3.2, and ~\cite{Los}, Theorem 6.10,  that 
$\Psi_{\mr{unif}}$
preserves the symplectic structure  up to a constant factor, i.e.,
\begin{align} \label{e06}
 \Psi_{\mr{unif}}^*(\omega_{\mfT^\Sigma}^{\mr{PGL}}) = \sqrt{-1} \cdot \omega^{\mr{can}}_{\mfT^\Sigma}.
 \end{align}
(Notice that we use  the conventions concerning $\omega_{\mfT^\Sigma}^{\mr{PGL}}$ chosen  in ~\cite{Los}. With these conventions, Kawai's result may be described as the equality (\ref{e06}), as B. Loustau mentioned in a footnote of {\it loc.\,cit.}.
In Kawai's original paper, he asserted the equality   $\Psi_{\mr{unif}}^*(\omega_{\mfT^\Sigma}^{\mr{PGL}}) = \pi \cdot \omega^{\mr{can}}_{\mfT^\Sigma}$.)

\vspace{3mm}
\subsection{} \label{s03}

Our aim in the present paper is to address the question
 {\it whether a similar result holds for hyperbolic (algebraic) curves of  positive characteristic}.
Just as in the case of the theory over $\mbC$, one may define the notion of an indigenous bundle in characteristic $p>0$ and their moduli space.
Various properties of such objects  were discussed in the context of the $p$-adic Teichm\"{u}ller  theory developed  by S. Mochizuki (cf. ~\cite{Mzk1}; ~\cite{Mzk2}).
One of the key ingredients in the development of this theory is
the study of  the $p$-curvature of indigenous bundles.
Recall that the $p$-curvature of a connection measures  the obstruction to the  compatibility of $p$-power structures that appear in certain associated spaces of infinitesimal (i.e., ``Lie'')  symmetries.
We say that an indigenous bundle is {\it dormant} (cf. Definition \ref{d01} (i)) if its $p$-curvature  vanishes identically.
This condition  implies that the underlying projective line bundle with connection is locally trivial in the Zariski  topology.

In many aspects (including the aspect just explained), dormant indigenous bundles  may be thought of as reasonable  (algebraic) products  used  to develop an analogous theory of  indigenous bundles on Riemann surfaces. 
As  explained in \S\,\ref{s02}, each connected  compact hyperbolic Riemann surface $X$ (with marking) of genus $g>1$ admits a canonical indigenous bundle $\mcP^\ind_X$ determined by  the section  $\sigma_{\mr{unif}}$.
Thus, the Teichm\"{u}ller space $\mfT^\Sigma$
   may be identified with  the moduli space classifying such $X$'s equipped with a specific  nice  indigenous bundle  (i.e., $\mcP^\ind_X$).
With that in mind,
we consider the moduli stack classifying proper hyperbolic curves of characteristic $p$ equipped with a dormant indigenous bundle  as a  characteristic $p$  analogue of such a covering space 
 of $\mfM_{g,\mbC}^\mr{an}$.
On the basis of this perspective, 
we give   an affirmative answer to the above question.

\vspace{3mm}
\subsection{} \label{s09}

In what follows, we shall  describe the main result of the present paper.
Let  $K$ be a field of characteristic $p>2$ and 
$g$ an integer $>1$.
Denote by 
\begin{align}
{^\circledcirc\Dp_{g, K}} \ (\text{resp.,} \ {^\circledcirc\mfS^{^\mr{Zzz...}}_{g, K}}) \notag
\end{align}
(cf. (\ref{e46}) and (\ref{e100}))   the moduli stack  classifying {\it ordinary dormant curves} (cf. Definition \ref{d01} (ii)) of genus $g$ over 
$K$ (resp., {\it ordinary dormant curves} of genus $g$ over $K$ equipped with an indigenous bundle).
Also, denote by
\begin{align}
{^\circledcirc T^{\vee ^\mr{Zzz...}}_{g,K}} \notag
\end{align}
(cf.  (\ref{e100}))
 the cotangent bundle 
 of  ${^\circledcirc\Dp_{g, K}}$.
 It is known (cf. Propositions \ref{p01} and  \ref{p05}) that   ${^\circledcirc\mfS^{^\mr{Zzz...}}_{g, K}}$ (resp., ${^\circledcirc T^{\vee ^\mr{Zzz...}}_{g,K}}$) may be represented by a  geometrically connected  smooth Deligne-Mumford stack over $K$ of dimension $6g-6$. 
As we will discuss in \S\,\ref{s25}, 
there exists  a canonical symplectic structure
 on ${^\circledcirc\mfS^{^\mr{Zzz...}}_{g, K}}$  (resp., ${^\circledcirc T^{\vee ^\mr{Zzz...}}_{g,K}}$), 
 which we denote by 
\begin{align}
\omega^{\mr{PGL}}_{\circledcirc} \ \ \  \left(\text{resp.}, \  \omega^\mr{can}_{\circledcirc} \right). \notag
\end{align}
Then, the  main result of the present paper is described as follows.

\vspace{3mm}
\begin{intthm}[cf. Theorem \ref{t05}] \label{tA}
The canonical  isomorphism
\begin{align}
\Psi_{g,K} : {^\circledcirc T^{\vee ^\mr{Zzz...}}_{g,K}} \isom {^\circledcirc\mfS^{^\mr{Zzz...}}_{g, K}}\notag
\end{align}
(cf. (\ref{e61}))  preserves  the  symplectic structure, i.e., 
\begin{align}
 \Psi_{g,K}^*(\omega^{\mr{PGL}}_\circledcirc) = \omega^\mr{can}_\circledcirc. \notag
  \end{align}
\end{intthm}
\vspace{3mm}

In particular,  the above theorem implies that 
the image of the canonical section 
${^\circledcirc\Dp_{g, K}} \migi  {^\circledcirc\mfS^{^\mr{Zzz...}}_{g, K}}$ is Langrangian (cf. Remark \ref{R071w}) with respect to the symplectic structure $\omega^{\mr{PGL}}_\circledcirc$.
Finally, as an application of Theorem A, 
we construct (cf. Corollary \ref{p76}) a so-called Frobenius-constant quantization on the moduli stack  
${^\circledcirc \mfS^{^\mr{Zzz...}}_{g,K}}$; such an additional structure  will be   of interest in the context of symplectic geometry.

\vspace{3mm}
\hspace{-4mm}{\bf Acknowledgements.}
 We express our  sincere and deepest gratitude to Professors Shinichi Mochizuki, Yuichiro Hoshi (and hyperbolic curves of positive characteristic!)
  for  their  helpful suggestions and heartfelt encouragement.
 Part of this work was written during the author's stay at Hausdorff Center for Mathematics in Bonn, Germany, whose
  hospitality and  very warm atmosphere is gratefully acknowledged.
The author's work  was partially  supported by  Grant-in-Aid for Scientific Research (KAKENHI) grant 24-5691 and by Grant-in-Aid for JSPS Fellows.
Finally,  thanks go to the referee for  a careful reading of the manuscript and for various useful comments and suggestions.
 \vspace{10mm}

\section{Preliminaries} \vspace{0mm}

In this section,  we shall review some definitions and facts concerning our discussion.


\vspace{3mm}
\subsection{Ground ring} \label{s305}

Throughout the present paper, we fix an integer $g>1$ and  a commutative ring $R$ over $\mbZ [\frac{1}{2}]$.
The assumption that $2$ is invertible in the ground ring $R$ will be necessary 
to construct  the symplectic structure $\omega^\mr{PGL}_{g, R}$ introduced later (cf. (\ref{e5f5})) and to  apply  previous   results on indigenous bundles in positive characteristic (cf. ~\cite{Mzk2}) to our discussion.
Also,  we shall write $\mfS \mfe \mft$ for the category of (small) sets.

\vspace{3mm}
\subsection{Sheaves and complexes} \label{dds305}

Let $S$ be a Deligne-Mumford stack over $R$   (cf. ~\cite{LMB} for the definition and basic properties of Deligne-Mumford stacks).
  Unless stated otherwise,   the structure sheaf $\mcO_S$ of  $S$ and
    all $\mcO_S$-modules  are considered as sheaves  in the \'{e}tale topology.
  If   $\nabla : \mcK^0 \migi \mcK^1$ is    a morphism of sheaves (in the \`{e}tale topology) of abelian groups on $S$, then
we often think of it as a complex concentrated at degree $0$ and $1$.
Denote this complex by
\begin{align} 
  \mcK^\bullet [\nabla].  \notag
  \end{align}
Also, any abelian sheaf $\mcF$ may be thought of as a complex concentrated at degree $0$.
For $n \in \mbZ$,  we define the complex
\begin{align} 
\mcF [n] \notag
\end{align}
to be $\mcF$ shifted down by $n$, so that $\mcF [n]^{-n} := \mcF$ and $\mcF [n]^i := 0$ for each $i \neq -n$.

If, moreover, $X$ is a  Deligne-Mumford stack over $S$,  
then we shall write $\Omega_{X/S}$ for the  sheaf  of 1-forms of $X$ over $S$,
 $\bigwedge^i\Omega_{X/S}$ ($i = 1,2, \cdots$) for its $i$-th exterior power, and 
 $\mcT_{X/S}$ for the dual $\mcO_X$-module $\Omega_{X/S}^\vee$  of $\Omega_{X/S}$ (i.e., the sheaf of derivations on  $\mcO_X$ over  $S$).

\vspace{3mm}
\subsection{Complex analytic spaces} \label{sd07}

Suppose that  $X$ is 
a complex analytic space  (resp., a complex analytic space over a complex analytic space $S$).
Then,  we shall write 
$\mcO_X$ for  the structure sheaf of $X$  consisting of holomorphic functions, 
$\Omega_X$ (resp., $\Omega_{X/S}$) for the sheaf of holomorphic $1$-forms of $X$ (resp., the sheaf of holomorphic $1$-forms  of $X$ over $S$), and 
  $\mcT_X$ (resp., $\mcT_{X/S}$) for the  dual $\mcO_X$-module of $\Omega_X$ (resp., $\Omega_{X/S}$).

Next, let us suppose that $X$ is  
  a scheme  (resp.,  a Deligne-Mumford stack) of finite type over $\mbC$.
  Then,  we shall write
$X^\mr{an}$ for  the complex analytic space (resp., the complex analytic stack, i.e., the  stack in groupoids over the category of complex analytic spaces equipped with the analytic topology) associated with $X$.

\vspace{3mm}
\subsection{Symplectic structures} \label{s07}

Let 
$X$ be a smooth Deligne-Mumford stack over $R$ of relative dimension $n >0$ (resp., a smooth  Deligne-Mumford  complex analytic stack (cf. ~\cite{BN}, Definition 3.3) of dimension $n >0$).
A {\bf symplectic structure}  on $X$  is, by definition,  a nondegenerate closed $2$-form $\omega \in \Gamma (X, \bigwedge^2\Omega_{X/R})$ (resp., $\omega \in \Gamma (X, \bigwedge^2 \Omega_X)$). 
Here,  a $2$-form $\omega$ is called  {\bf nondegenerate}  if the morphism $\Omega_{X/R} \migi \mcT_{X/R}$ (resp., $\Omega_{X} \migi \mcT_X$) induced naturally by $\omega$ is an isomorphism.

 Given an $\mcO_X$-module $\mcF$,
we shall write
\begin{align} 
\mbA (\mcF) \notag
\end{align}
for the total space of $\mcF$.
In the non-resp'd case, $\mbA (\mcF)$ forms the relative affine scheme $\mcS pec (\mbS (\mcF^\vee))$ over $X$,  where $\mbS(\mcF^\vee)$ denotes the symmetric algebra on the dual $\mcF^\vee$  of $\mcF$ over $\mcO_X$.
Denote by $T^\vee_X$ the total space of  $\Omega_{X/R}$ (resp., $\Omega_X$), i.e., 
\begin{align}
   T^\vee_X := \mbA (\Omega_{X/R})  \  \left(\text{resp.,} \  T^\vee_X:= \mbA (\Omega_{X}) \right), \notag
   \end{align}
which is a smooth Deligne-Mumford stack over $R$ of relative dimension $2 n$ (resp., a smooth  Deligne-Mumford  complex analytic stack of dimension $2n$);
we shall refer to $T^\vee_X$ as the {\bf cotangent bundle} of $X$.
Denote by
\begin{align} 
\pi_{X}^{T^\vee} : T^\vee_X \migi X, \hspace{5mm} 0_X : X \migi T^\vee_X \notag
\end{align}
the natural projection and the zero section respectively.
It is well-known that there exists a unique $1$-form 
\begin{align} 
\lambda_X \in \Gamma (T^\vee_X, \Omega_{T^\vee_X /R}) \  \ (\text{resp.,} \  \lambda_X \in \Gamma (T^\vee_X, \Omega_{T^\vee_X})) \notag
\end{align}
 on   $T^\vee_X$ determined by  the following condition:
if $\lambda_\sigma$ is the $1$-form on an open subscheme $U$ of $X$ corresponding to a local section $\sigma : U \migi T^\vee_X$ of $\pi_X^{T^\vee}$, then the equality   $\sigma^*(\lambda_X) = \lambda_\sigma$ holds.
We shall refer to $\lambda_X$ as the {\bf Liouville form} on $T^\vee_X$.
By construction, the Liouville form $\lambda_X$ lies in $\Gamma (T^\vee_X, \pi^{T^\vee *}_X (\Omega_{X/R})) \subseteq \Gamma (T^\vee_X, \Omega_{T^\vee_X/R})$
(resp., $\Gamma (T^\vee_X, \pi^{T^\vee *}_X (\Omega_{X})) \subseteq \Gamma (T^\vee_X, \Omega_{T^\vee_X})$).
Its exterior derivative 
\begin{align} \label{e10}
   \omega^{\mr{can}}_{X} := d\lambda_X \in \Gamma (T^\vee_X, {\bigwedge}^2\Omega_{T^\vee_X/R}) \ \left(\text{resp.,} \   \omega^{\mr{can}}_{X} := d\lambda_X\in \Gamma (T^\vee_X, {\bigwedge}^2\Omega_{T^\vee_X})\right)
   \end{align}
defines a symplectic structure on $T^\vee_X$.
If $q_1, \cdots, q_n$ are local coordinates in $X$,
  then the dual coordinates $p_1, \cdots,  p_n$ in $T^\vee_X$ are the coefficients of the decomposition of the $1$-form $\lambda_X$ into linear combination of the differentials $dq_i$, i.e.,
$\lambda_X = \sum_{i=1}^n p_idq_i$.
Hence, $\omega^{\mr{can}}_{X}$ may be expressed locally as
$\omega^{\mr{can}}_{X} = \sum_{i=1}^n dp_i \wedge dq_i$.

\vspace{3mm}

\subsection{Curves and their moduli space} \label{s08}

Let  $S$ be  a Deligne-Mumford  stack   over $R$.
By a {\bf curve} over $S$, we mean a geometrically connected    smooth relative   scheme $f : X \migi S$ over $S$ of relative dimension $1$.
Here, a morphism $f : X \migi S$ between Deligne-Mumford stacks  is called a {\it relative scheme} if it is a schematic morphism.
Also, we shall say that a proper  curve $f :X \migi S$ over $S$ is  {\bf of genus $g$} 
if the direct image $f_*(\Omega_{X/S})$ of $\Omega_{X/S}$  is a locally free $\mcO_S$-module of constant rank $g$.
 Write 
\begin{align} 
\mfM_{g, R} \notag
\end{align}
 for the moduli stack classifying  proper curves of genus $g$  over $R$; it may be represented by  a geometrically connected smooth Deligne-Mumford stack over $R$ of relative dimension $3g-3$.
 Also,  write 
\begin{align} 
f_{g, R} : C_{g,R} \migi\mfM_{g,R} \notag
\end{align}
 for the tautological curve over $\mfM_{g,R}$ and
 \begin{align}
 \mfS \mfc \mfh_{/ \mfM_{g,R}} \notag
 \end{align}
for the category of relative schemes over $\mfM_{g,R}$.

It follows from Serre duality that for any proper curve $f : X \migi S$,  the $\mcO_S$-module $\mbR^1f_{*}(\Omega_{X/S})$ is isomorphic to $\mcO_{S}$.
Throughout the present paper, we fix a specific choice of an isomorphism
\begin{align} 
   \int_{C_{g,\mbZ[\frac{1}{2}]}} : \mbR^1f_{g,\mbZ[\frac{1}{2}]*}(\Omega_{C_{g,\mbZ[\frac{1}{2}]}/\mfM_{g,\mbZ[\frac{1}{2}]}})\isom \mcO_{\mfM_{g,\mbZ[\frac{1}{2}]}}   \notag
   \end{align}
 of $\mcO_{\mfM_{g,\mbZ[\frac{1}{2}]}}$-modules (i.e., the {\it trace map}).
 Note that the isomorphism $ \int_{C_{g,\mbZ[\frac{1}{2}]}}$ will be used   to determine both  the symplectic structure $\omega^\mr{PGL}_{g, R}$ (cf. (\ref{e5f5})) and the morphism denoted by  $\Psi_{g, K}$ (cf. (\ref{e61});  we have to choose them in such a coherent way in order to obtain the equality asserted in Theorem A.

If $u : T \migi \mfM_{g,R}$ is an object of  $\mfS \mfc \mfh_{/ \mfM_{g,R}}$, 
 then we shall write 
\begin{align} 
f_T : C_T \migi T \notag
\end{align}
 for the curve over $T$ classified by $u$, i.e., $C_T := C_{g,R} \times_{f_{g,R}, \mfM_{g,R}, u} T$.
 We obtain an isomorphism
 \begin{align} 
  \int_{C_T} : \mbR^1f_{T*}(\Omega_{C_T/T})\isom \mcO_T  \notag
  \end{align}
defined as   the pull-back of  $ \int_{C_{g,\mbZ[\frac{1}{2}]}}$ via the composite $T \xrightarrow{u} \mfM_{g,R} \migi \mfM_{g,\mbZ[\frac{1}{2}]}$.
For a vector bundle $\mcE$ on $C_T$ (i.e., a locally free coherent $\mcO_{C_T}$-module),  
denote by 
 \begin{align} \label{e16}  
   \int_{C_T, \mcE} : \mbR^1f_{T*}(\Omega_{C_T/T} \otimes \mcE^\vee) \isom f_{T*}(\mcE)^\vee   
 \end{align}
the isomorphism induced  from  the pairing 
\begin{align} 
\mbR^1f_{T*}( \Omega_{C_T/T} \otimes \mcE^\vee) \otimes f_{T*}(\mcE)  &  \xrightarrow{\hspace{8mm} \cup \hspace{8mm}} 
\mbR^1f_{T*}(\Omega_{C_T/T} \otimes (\mcE^\vee \otimes \mcE)) \notag \\
& \hspace{-1mm} \xrightarrow{\hspace{19mm}}
 \mbR^1 f_{T*} (\Omega_{C_T/T}) \notag  \\
& \hspace{-0mm}\xrightarrow{\hspace{6mm} \int_{C_T} \hspace{6mm}} \mcO_T, \notag
\end{align}
where the second arrow denotes the morphism arising from the natural pairing $\mcE^\vee \otimes \mcE \migi \mcO_{C_T}$.

\vspace{3mm}
\subsection{Connections} \label{s09}
Let 
$S$ be a relative scheme over $\mfM_{g, R}$, i.e., an object of $\mfS \mfc \mfh_{/\mfM_{g,R}}$,
 which classifies  
a proper curve $f_S : C_S \migi S$ over $S$.
Let $G$  be a connected  smooth algebraic group over $R$ with Lie algebra $\mfg$ and 
  $\pi : \mcP \migi C_S$  a (right) $G$-torsor over $C_S$.
 Write  
 \begin{align}
 \mr{ad}(\mcP) 
 \end{align}
   for  the adjoint vector bundle associated with  $\mcP$.
   That is to say, 
   $\mr{ad}(\mcP)$ is the vector bundle obtained from $\mcP$ via change of structure group by 
    the adjoint representation $\mr{Ad} : G \migi \mr{GL} (\mfg)$.
    The direct image $\pi_*(\mcT_{\mcP/S})$ has a $G$-action arising from the $G$-action on $\mcP$, and hence,  we obtain  its subsheaf
   \begin{align}
   \widetilde{\mcT}_{\mcP/S} \ \left( := (\pi_*(\mcT_{\mcP/S}))^G \right)\notag
   \end{align}
    consisting   of $G$-invariant sections.
The differential of  $\pi$ gives  a short exact sequence of $\mcO_{C_S}$-modules
\begin{align} \label{e17}
  0 \migi \mr{ad}(\mcP) \migi  \widetilde{\mcT}_{\mcP/S} \xrightarrow{\alpha_\mcP} \mcT_{C_S/S} \migi 0.
  \end{align}
An {\bf $S$-connection} on $\mcP$ is, by definition,  a split injection $\nabla_\mcP : \mcT_{C_S/S} \migi  \widetilde{\mcT}_{\mcP/S}$ of the  short exact sequence (\ref{e17}), i.e., $\alpha_\mcP \circ \nabla_\mcP = \mr{id}$. 
Since  $C_S$ is of relative dimension $1$ over $S$,  any such $S$-connection is necessarily {\it integrable},  which means that it is 
compatible with the respective Lie bracket structures on $\mcT_{C_S/S}$ and $\widetilde{\mcT}_{\mcP/S}= (\pi_*(\mcT_{\mcP/S}))^G$.
By a {\bf flat $G$-torsor}  over $C_S$, we mean a pair $(\mcP, \nabla_\mcP)$ consisting of a $G$-torsor $\mcP$ over $C_S$ and an $S$-connection $\nabla_\mcP$ on $\mcP$.

If $G = \mr{GL}_m$ for some $m \geq 1$, 
then the notion of an $S$-connection on a $\mr{GL}_m$-torsor $\mcP$ recalled  here may be identified with
the usual definition of an $S$-connection  (cf. ~\cite{Kal}, \S\,1.0) on the associated vector bundle $\mcP \times^{\mr{GL}_m} (\mcO_{C_S}^{\oplus m})$ (cf. Remark \ref{R071} for a detailed discussion);
 in this situation, we shall not distinguish between these notions of connections.

Given  an $S$-connection $\nabla_\mcP$ on $\mcP$,  we denote by 
\begin{align} \label{e18}
   \nabla^\mr{ad}_\mcP : \mr{ad}(\mcP) \migi \Omega_{C_S/S} \otimes \mr{ad}(\mcP) 
   \end{align}
the $S$-connection on $\mr{ad}(\mcP)$ induced  by $\nabla_\mcP$ via the  change of structure group by $\mr{Ad} : G \migi \mr{GL} (\mfg)$.
More explicitly, $\nabla_\mcP^{\mr{ad}}$ is the connection uniquely determined by the condition that 
$\langle \partial_1, \nabla_\mcP^{\mr{ad}} (\partial_2) \rangle = [\nabla_\mcP (\partial_1), \partial_2]$ for any local sections $\partial_1 \in \mcT_{C_S/S}$ and  $\partial_2 \in \mr{ad} (\mcP)$, where $\langle -, - \rangle$ denotes the pairing $\mcT_{C_S/S} \times (\Omega_{C_S/S} \otimes \mr{ad} (\mcP)) \migi \mr{ad} (\mcP)$ arising from the natural pairing $\mcT_{C_S/S} \times \Omega_{C_S/S} \migi \mcO_{C_S}$.

\vspace{3mm}
\subsection{Ring of differential operators} \label{sss13}

Recall (cf.  ~\cite{Be3}, \S\,1.2) that  the sheaf of {\it   crystalline differential operators} 
 on $C_S$ over $S$ is the Zariski sheaf 
\begin{equation} \mcD_{C_S/S}\notag\end{equation}
 on $C_S$
   generated, as a sheaf of noncommutative   rings, by $\mcO_{C_S}$ and  
$\mcT_{C_S/S}$
 subject to the  relations
\begin{align}
f_1 \ast f_2 = f_1 \cdot f_2,  \hspace{6mm}
f_1 \ast \xi_1 = f_1 \cdot \xi_1, \hspace{6mm}
\xi_1 \ast \xi_2 -\xi_2 \ast \xi_1=  [\xi_1, \xi_2], \hspace{6mm}
\xi_1 \ast f_1 - f_1 \ast \xi_1 = \xi_1 (f_1) \notag
\end{align}
for any local sections $f_1$, $f_2 \in \mcO_{C_S}$ and $\xi_1$, $\xi_2 \in \mcT_{C_S/S}$, where $\ast$ denotes the multiplication in $ \mcD_{C_S/S}$.
In a usual sense, the {\it order $(\geq 0)$} of a given  crystalline differential operator, i.e., a 
 local section of $\mcD_{C_S/S}$, 
 is well-defined.
Hence, $\mcD_{C_S/S}$ admits, for each $j \geq 0$, 
 the subsheaf
 \begin{equation}
 \mcD_{C_S/S}^{\leq  j} \ \left(\subseteq  \mcD_{C_S/S}\right)\notag
 \end{equation}  consisting of local sections of $ \mcD_{C_S/S}$  of order $\leq  j$.
The sheaf $\mcD_{C_S/S}$ (resp.,  $\mcD^{< j}_{C_S/S}$ for each $j = 0, 1, 2, \cdots$)   admits  two different  structures  of  $\mcO_{C_S}$-module --- one as given by left multiplication, where we denote this $\mcO_{C_S}$-module by $^l\mcD_{C_S/S}$ (resp., $^l\mcD^{\leq   j}_{C_S/S}$), and  the other given by right multiplication, where  we denote this $\mcO_{C_S}$-module by $^r\mcD_{C_S/S}$ (resp., $^r\mcD^{\leq j}_{C_S/S}$) ---.
In particular,  we have ${^l\mcD}_{C_S/S}^{\leq 0} = {^r\mcD}_{C_S/S}^{\leq 0} = \mcO_{C_S}$ (as $\mcO_{C_S}$-modules).
The set  $\{ \mcD_{C_S/S}^{\leq  j} \}_{j\geq 0}$ forms an increasing filtration on $\mcD_{C_S/S}$ satisfying that
 \begin{equation} \label{GL7}
 \bigcup_{j \geq 0} \mcD_{C_S/S}^{\leq  j} = \mcD_{C_S/S},  \ \text{and} \  \ \mcD_{C_S/S}^{\leq j}/\mcD_{C_S/S}^{\leq (j-1)} \cong \mcT_{C_S/S}^{\otimes j} \ \text{ for every $j \geq 1$.}
 \end{equation} 

Let $\mcF$ be a  vector bundle on $C_S$.
In what follows,  
we shall regard  the tensor product  ${\mcD}^{\leq j}_{C_S/S} \otimes \mcF:= {^r\mcD}^{\leq j}_{C_S/S} \otimes \mcF$ (resp., $\mcF \otimes \mcD^{\leq j}_{C_S/S} := \mcF \otimes {^l \mcD}^{\leq j}_{C_S/S}$)
 as being equipped with a  structure of  $\mcO_{C_S}$-module arising from the structure of $\mcO_{C_S}$-module ${^l\mcD}^{\leq j}_{C_S/S}$ (resp., ${^r \mcD}^{\leq j}_{C_S/S}$)
   on ${\mcD}^{\leq j}_{C_S/S}$.

Next, let $\nabla_\mcF$ be an $S$-connection on $\mcF$.
 It induces   a structure of left ${\mcD}_{C_S/S}$-module 
 \begin{equation}
 \widehat{\nabla}_\mcF : {\mcD}_{C_S/S} \otimes \mcF \migi \mcF\notag
 \end{equation}
  on $\mcF$  determined uniquely by the condition that
$\widehat{\nabla}_\mcF (\partial \otimes v) = \langle  \partial,   \nabla_\mcF (v)\rangle$ for any local sections $v \in \mcF$ and $\partial \in \mcT_{C_S/S}$, where $\langle -,- \rangle$ denotes the
pairing $\mcT_{C_S/S} \times (\Omega_{C_S/S}\otimes \mcF) \migi \mcF$ induced by
the natural paring $\mcT_{C_S/S} \times \Omega_{C_S/S} \migi \mcO_{C_S}$.
The assignment $\nabla_\mcF \mapsto \widehat{\nabla}_\mcF$ determines 
 a  bijective correspondence between the set of $S$-connections 
  on $\mcF$ and the set of 
  structures of left ${\mcD}_{C_S/S}$-module 
  on $\mcF$.

\vspace{3mm}
\begin{rema} \label{R071}
 As mentioned in \S\,\ref{s09}, there exists 
  a bijective correspondence between   $S$-connections
 on a $\mr{GL}_m$-torsor $\mcP$ and  $S$-connections   (in the classical sense) on the corresponding vector bundle $\mcF := \mcP \times^{\mr{GL}_m} (\mcO_{C_S}^{\oplus m})$.
 In this remark, we describe this correspondence in somewhat detail.
 
 Let us write
   \begin{align}
  \mcD {\it iff}^{\leq 1}_{\mcF, \mcF} := \mcH om_{\mcO_{C_S}} (\mcF, \mcF \otimes \mcD_{C_S/S}^{\leq 1}), \notag
  \end{align}
i.e., the sheaf of first order differential operators from $\mcF$ to $\mcF$ itself.
  The natural surjection  $\mcD^{\leq 1}_{C_S/S} \migisurj \mcT_{C_S/S}$ (cf. (\ref{GL7})) induces
  a surjection 
 \begin{align}
 q : \mcD {\it iff}^{\leq 1}_{\mcF, \mcF}  \migisurj \mcE nd_{\mcO_{C_S}} (\mcF) \otimes \mcT_{C_S/S}.\notag
 \end{align}
Note that the morphism   $\mcT_{C_S/S} \migiincl \mcE nd_{\mcO_{C_S}} (\mcF) \otimes \mcT_{C_S/S}$ given by assigning $\partial \mapsto  \mr{id}_\mcF \otimes \partial$ for any local section $\partial \in \mcT_{C_S/S}$ is injective;
 we shall regard  
  $\mcT_{C_S/S}$ as an $\mcO_{C_S}$-submodule of $\mcE nd_{\mcO_{C_S}} (\mcF) \otimes \mcT_{C_S/S}$ by this injection.
 Let us write
  \begin{align}
  \mcD {\it iff}^{\leq 1 \blacklozenge}_{\mcF, \mcF} := q^{-1} (\mcT_{C_S/S}) \ \left(\subseteq  \mcD {\it iff}^{\leq 1}_{\mcF, \mcF} \right), \notag
    \end{align}
  which admits a surjection
\begin{align}
\alpha_\mcF \ \left(:= q |_{ \mcD {\it iff}^{\leq 1 \blacklozenge}_{\mcF, \mcF}}\right) 
: \mcD {\it iff}^{\leq 1 \blacklozenge}_{\mcF, \mcF} \migisurj \mcT_{C_S/S}.\notag
\end{align}
 Recall  (cf. ~\cite{BiSu}, \S\,2, (2.2))  that there exists a canonical isomorphism
  \begin{align}
  \xi  :  \widetilde{\mcT}_{\mcP/S} \isom \mcD {\it iff}^{\leq 1 \blacklozenge}_{\mcF, \mcF} \notag
  \end{align}
of $\mcO_{C_S}$-modules making the following diagram commute:
\begin{align}
\xymatrix{
\widetilde{\mcT}_{\mcP/S} \ar[rr]^-{\xi}_-{\sim} \ar[rd]_-{\alpha_\mcP}& & \mcD {\it iff}^{\leq 1 \blacklozenge}_{\mcF, \mcF}  \ar[ld]^-{\alpha_\mcF}
\\
& \mcT_{C_S/S}.&
}\notag
\end{align}

Now, suppose that we are given an $S$-connection $\nabla_\mcP^1 : \mcT_{C_S/S} \migi \widetilde{\mcT}_{\mcP/S}$ on the $\mr{GL}_m$-torsor  $\mcP$.
Then, one may find an $S$-connection $\nabla^1_\mcF : \mcF \migi \Omega_{C_S/S} \otimes \mcF$ on $\mcF$ satisfying the equality
$(\xi \circ \nabla^1_\mcP) (\partial) (v) - v \otimes \partial=  \langle \partial, \nabla^1_\mcF (v) \rangle$
for any local sections $\partial \in \mcT_{C_S/S}$ and $v \in \mcF$, where $\langle -, - \rangle$ denotes the pairing $\mcT_{C_S/S} \times (\Omega_{C_S/S} \otimes \mcF) \migi \mcF$ arising from the natural pairing $\mcT_{C_S/S} \times \Omega_{C_S/S} \migi \mcO_{C_S}$.

Conversely, let $\nabla_\mcF^2 : \mcF \migi \Omega_{C_S/S} \otimes \mcF$ be an $S$-connection on the vector bundle  $\mcF$ and  $\partial$ a local section of  $\mcT_{C_S/S}$.
Then, the assignment $v \mapsto  \langle \partial,  \nabla^2_\mcF (v) \rangle + v \otimes \partial$ (where  $v \in \mcF$) determines a local section $\delta_\partial$  of $\mcD {\it iff}^{\leq 1 \blacklozenge}_{\mcF, \mcF}$.
The morphism $\nabla^2_\mcP : \mcT_{C_S/S} \migi \widetilde{\mcT}_{\mcP/S}$ given by assigning $\partial \mapsto \xi^{-1} (\delta_\partial)$ forms an $S$-connection on $\mcP$.

One verifies that the assignments $\nabla^1_\mcP \mapsto \nabla^1_\mcF$ and $\nabla^2_\mcF \mapsto \nabla^2_\mcP$ determine a bijective correspondence between the set of $S$-connections on $\mcP$ and the set of $S$-connections on $\mcF$, as  desired.
\end{rema}

\vspace{10mm}
\section{Indigenous Bundles} \label{s10}\vspace{3mm}

In this section, we recall the notion of an indigenous bundle on a curve
 and some properties concerning  indigenous bundles.
We refer the reader to ~\cite{Mzk1}, ~\cite{Mzk2},  ~\cite{Wak}, ~\cite{Wak5}
for detailed discussions involved.

\vspace{3mm}
\subsection{} \label{s11}

Let us fix 
a  relative scheme $S$ over $\mfM_{g, R}$,
 which classifies 
 a proper curve $f_S : C_S \migi S$ of genus $g$. 
In what follows, let   $G$  be  the projective linear group over $R$ of rank $2$, i.e.,  $G := \mr{PGL}_2$.
Since $2$ is invertible in $R$, its Lie algebra $\mfg$ is naturally isomorphic to  $\mfs \mfl_2$.
Write $B$ for the  Borel subgroup of $G$ defined to be the image  of upper triangular 
 matrices via the natural  quotient $\mr{GL}_2 \migisurj G$.
We recall from
   ~\cite{F}, \S4, or ~\cite{Mzk1}, Chap.\,I, \S\,2, Definition 2.2,
 the following definition of an indigenous bundle: 

\vspace{3mm}
\bde  \label{D01}
\begin{itemize}
\item[(i)]
Let $\mcP^\ind :=(\mcP_B, \nabla_{\mcP_G})$ be a  pair consisting of a  (right) $B$-torsor $\mcP_B$ over $C_S$ and  an  $S$-connection $\nabla_{\mcP_G}$ on the $G$-torsor $\mcP_{G} := \mcP_B \times^B G$ induced by $\mcP_B$.
We shall say that $\mcP^\ind$ is an  {\bf indigenous bundle} on $C_S/S$ if 
the composite
\begin{align} \label{e21}
 \overline{\nabla}_{\mcP_G}: \mcT_{C_S/S} \xrightarrow{\nabla_{\mcP_G}}\widetilde{\mcT}_{\mcP_{G}/S} \migisurj  \widetilde{\mcT}_{\mcP_{G}/S}  / \widetilde{\iota}(\widetilde{\mcT}_{\mcP_{B}/S})
 \end{align}
is an isomorphism, where $\widetilde{\iota}$ denotes the natural injection $\widetilde{\mcT}_{\mcP_{B}/S} \migiincl  \widetilde{\mcT}_{\mcP_{G}/S}$.
\item[(ii)]
Let $\mcP^\ind := (\mcP_B, \nabla_{\mcP_G})$ and $\mcQ^\ind :=(\mcQ_B, \nabla_{\mcQ_G})$ be  indigenous bundles on $C_S/S$.
 An {\bf isomorphism of indigenous bundles} from  $\mcP^\ind$ to $\mcQ^\ind$ is an isomorphism $\mcP_B \isom \mcQ_B$ of $B$-torsors  such that
 the induced isomorphism $\mcP_\G \isom \mcQ_\G$ of $\G$-torsors 
  is compatible with the respective $S$-connections.
  \end{itemize}
  \ede
\vspace{3mm}

Then,  we already know that the following proposition holds.

\vspace{3mm}
\bpr [cf. ~\cite{Mzk1}, Chap.\,I, \S\,2, Theorem 2.8] \label{R002} 
Any indigenous bundle on $C_S/S$ does not have  nontrivial automorphisms.
 \epr

\vspace{3mm}
\subsection{} \label{s12}

In what follows, we shall  construct a canonical  filtration on the adjoint vector  bundle associated with  the underlying $G$-torsor of  an indigenous bundle.
Let  $\mcP^\ind := (\mcP_B, \nabla_{\mcP_G})$ be an indigenous bundle on $C_S/S$.
Consider the morphism of short exact sequences
\begin{align} 
  \begin{CD} 0 @>>>\mr{ad}(\mcP_B) @>>> \widetilde{\mcT}_{\mcP_B/S} @> \alpha_{\mcP_B} >> \mcT_{C_S/S} @>>>0
\\
@. @VV \iota V @VV \widetilde{\iota} V @VV \mr{id} V @.
\\
0 @>>>  \mr{ad}(\mcP_\G) @>>> \widetilde{\mcT}_{\mcP_\G/S} @>>\alpha_{\mcP_G} > \mcT_{C_S/S} @>>>0 \end{CD}\notag\end{align}
arising from the change  of  structure group $B \migiincl  G$.
This diagram yields an isomorphism
\begin{align} \label{e23}
 \mr{ad}(\mcP_\G) / \iota(\mr{ad}(\mcP_B)) \isom   \widetilde{\mcT}_{\mcP_\G/S} /\widetilde{\iota}(\widetilde{\mcT}_{\mcP_B/S}).
   \end{align}

Let us   define a $3$-step decreasing filtration $\{ \mr{ad}(\mcP_\G)^i\}_{i=0}^3$ on the rank $3$ vector bundle $\mr{ad}(\mcP_G)$ as follows:
\begin{align} 
  \mr{ad}(\mcP_G)^0 := & \   \mr{ad}(\mcP_G), \notag \\
 \mr{ad}(\mcP_G)^1:= & \ \iota(\mr{ad}(\mcP_B)), \notag \\
 \mr{ad}(\mcP_G)^2 := & \ \mr{Ker}\left( \mr{ad}(\mcP_G)^1 \xrightarrow{\nabla_{\mcP_G}^\mr{ad}|_{\mr{ad}(\mcP_G)^1}}  \Omega_{C_S/S}\otimes \mr{ad}(\mcP_G) 
  \migisurj  \Omega_{C_S/S}\otimes (\mr{ad}(\mcP_G)/ \mr{ad}(\mcP_G)^1) \right), \notag \\
  \mr{ad}(\mcP_G)^3 := & \ 0 \notag
 \end{align}
(cf. (\ref{e18}) for the definition of $\nabla_{\mcP_G}^\mr{ad}$).
It follows from the definition of an indigenous bundle that 
 \begin{align} 
     \nabla_{\mcP_G}^\mr{ad}(\mr{ad}(\mcP_G)^{j+1}) \subseteq \Omega_{C_S/S} \otimes \mr{ad}(\mcP_G)^j,  \notag\end{align}
for any $j \in \{0, 1\}$
and  the  {\it $\mcO_{C_S}$-linear} morphism
\begin{align} 
  \overline{\nabla}^{\mr{ad}, j+1}_{\mcP_G}  : \mr{ad}(\mcP_G)^{j+1} / \mr{ad}(\mcP_G)^{j+2} \migi \Omega_{C_S/S}\otimes (\mr{ad}(\mcP_G)^{j}/\mr{ad}(\mcP_G)^{j+1} ) \notag
  \end{align}
induced by  $\nabla_{\mcP_G}^\mr{ad}$
  is an isomorphism.
Denote by
\begin{align} 
    \overline{\nabla}^{\flat}_{\mcP_G} :\left(\mr{ad}(\mcP_\G) /\mr{ad}(\mcP_\G)^1 = \right) \  \mr{ad}(\mcP_\G) / \iota(\mr{ad}(\mcP_B))   \isom \mcT_{C_S/S}  \notag   \end{align}
the composite of  (\ref{e23}) and  $\overline{\nabla}_{\mcP_G}^{-1} :   \widetilde{\mcT}_{\mcP_\G/S} /\widetilde{\iota}(\widetilde{\mcT}_{\mcP_B/S}) \isom \mcT_{C_S/S} $ (cf. (\ref{e21})).
Also, denote by $\overline{\nabla}^\sharp_{\mcP_G}$ the composite isomorphism 
\begin{align}
\overline{\nabla}^\sharp_{\mcP_G} : \Omega_{C_S/S}  &\isom \Omega_{C_S/S} \notag\\
&\hspace{-0.5mm} \xrightarrow{\mr{id}_{\Omega_{C_S/S}^{\otimes 2}}\otimes (\overline{\nabla}^{\flat}_{\mcP_G})^{-1}} \Omega_{C_S/S}^{\otimes 2} \otimes ( \mr{ad}(\mcP_\G)/ \mr{ad}(\mcP_\G)^1) \notag  \\
&  \hspace{-0.5mm} \xrightarrow{\mr{id}_{\Omega_{C_S/S}}\otimes (\overline{\nabla}_{\mcP_G}^{\mr{ad}, 1})^{-1}} \Omega_{C_S/K} \otimes ( \mr{ad}(\mcP_\G)^1/ \mr{ad}(\mcP_\G)^2) \notag  \\
&  \hspace{-0.5mm}  \xrightarrow{(\overline{\nabla}_{\mcP_G}^{\mr{ad}, 2})^{-1}}  \mr{ad}(\mcP_\G)^2, \notag
\end{align}
where the first arrow denotes the automorphism of $\Omega_{C_S/S}$ given by multiplication by $2$, which is invertible in $R$ by assumption.
Note that  the first arrow in this composite should  be  added in order
to conclude  Lemma \ref{p40} described later.

\vspace{3mm}
\subsection{} \label{s15}

Denote by 
\begin{align}
\widetilde{\nabla}^\mr{ad}_{\mcP_G} : \widetilde{\mcT}_{\mcP_G/S} \migi \Omega_{C_S/S} \otimes \mr{ad} (\mcP_G)  \notag
\end{align}
a unique  
$f_S^{-1}(\mcO_S)$-linear 
morphism determined  by the condition that
\begin{align}
\langle \partial_1, \widetilde{\nabla}_{\mcP_G}^\mr{ad} (\partial_2)  \rangle = [\nabla_{\mcP_G} (\partial_1), \partial_2] - \nabla_\mcE ([\partial_1, \alpha_{\mcP_G} (\partial_2)])  \notag
\end{align}
 for any local sections $\partial_1 \in \mcT_{C_S/S}$ and $\partial_2 \in \widetilde{\mcT}_{\mcP_G/S}$,
  where $\langle -, - \rangle$ denotes the pairing $\mcT_{C_S/S} \times (\Omega_{C_S/S} \otimes \mr{ad} (\mcP_G)) \migi \mr{ad} (\mcP_G)$ arising from the natural pairing  $\mcT_{C_S/S} \times \Omega_{C_S/S} \migi \mcO_{C_S}$.
This morphism fits into the following short exact sequence of complexes:
\begin{align} \label{Ww500}
\begin{CD}
0 @>>> \mr{ad}(\mcP_G) @>>> \widetilde{\mcT}_{\mcP_G/S}@> \alpha_{\mcP_G}>> \mcT_{C_S/S}@>>> 0
\\
@. @VV \nabla^\mr{ad}_{\mcP_G}V @VV \widetilde{\nabla}^\mr{ad}_{\mcP_G} V @VVV @.
\\
0 @>>> \Omega_{C_S/S}\otimes\mr{ad}(\mcP_G)
@>>\mr{id} > \Omega_{C_S/S}\otimes \mr{ad}(\mcP_G)  @>  >> 0 @>>> 0,
\end{CD}
\end{align}
where the upper horizontal sequence is (\ref{e17}). 
Since  $\widetilde{\nabla}^\mr{ad}_{\mcP_G} \circ \nabla_{\mcP_G} =0$,
the connection $\nabla_{\mcP_G} : \mcT_{C_S/S} \migi \widetilde{\mcT}_{\mcP_G/S}$ induces a morphism $\mcT_{C_S/S}[0] \migi \mcK^\bullet [\widetilde{\nabla}^\mr{ad}_{\mcP_G}]$ defining  a split injection of (\ref{Ww500}).
This implies that
the short exact sequence
\begin{align} \label{Ww200}
0 \longmigi \mbR^1 f_{S*}(\mcK^\bullet [\nabla^\mr{ad}_{\mcP_G}])
\longmigi 
 \mbR^1 f_{S*}(\mcK^\bullet [\widetilde{\nabla}^\mr{ad}_{\mcP_G}])
\longmigi \mbR^1 f_{S*}(\mcT_{C_S/S}) \longmigi 0
\end{align}
obtained from (\ref{Ww500}) by applying the functor $\mbR^1 f_{S*}(-)$
is exact.

Next, let us consider 
the restriction
\begin{align} \label{Ww505}
\widetilde{\nabla}_{\mcP_B}^\mr{ad} \ \left(:= \widetilde{\nabla}_{\mcP_G}^\mr{ad}|_{\widetilde{\mcT}_{\mcP_B}} \right) : \widetilde{\mcT}_{\mcP_B} \migi\Omega_{C_S/S} \otimes \mr{ad} (\mcP_G)
\end{align}
of $\widetilde{\nabla}^\mr{ad}_{\mcP_G}$;
it fits into  the following morphism of short exact sequences:
\begin{align} \label{W51}
\begin{CD}
0 @>>> \widetilde{\mcT}_{\mcP_B/S} @>\mr{incl.}>> \widetilde{\mcT}_{\mcP_G/S}  @> \alpha'_{\mcP_G}   >> \mcT_{C_S/S} @>>> 0
\\
@. @VV \widetilde{\nabla}^\mr{ad}_{\mcP_B} V @VV \widetilde{\nabla}^\mr{ad}_{\mcP_G} V @VVV @.
\\
0 @>>> \Omega_{C_S/S}\otimes\mr{ad}(\mcP_G)@>> \mr{id} > \Omega_{C_S/S} \otimes \mr{ad}(\mcP_G)  @>>> 0 @>>> 0,
\end{CD}
\end{align}
where 
$\alpha'_{\mcP_G}$
 denotes the composite of the natural  quotient $\widetilde{\mcT}_{\mcP_G/S} \migisurj \widetilde{\mcT}_{\mcP_G/S}/\widetilde{\iota}(\widetilde{\mcT}_{\mcP_B/S})$ 
 and 
  $\overline{\nabla}_{\mcP_G}^{-1} :  \widetilde{\mcT}_{\mcP_G/S}/\widetilde{\iota}(\widetilde{\mcT}_{\mcP_B/S}) \isom \mcT_{C_S/S}$ (cf. (\ref{e21})).
Since $f_{S*}(\mcT_{C_S/S}) =0$, 
  the diagram (\ref{W51}) induces 
  an injection
  \begin{align}\label{W7000}
   \mbR^1 f_{S*} (\mcK^\bullet [\widetilde{\nabla}^\mr{ad}_{\mcP_B}]) \stackrel{}{\migiincl} 
\mbR^1f_{S*} (\mcK^\bullet [\widetilde{\nabla}^\mr{ad}_{\mcP_G}]).
   \end{align}

\vspace{3mm}
\subsection{} \label{s154}

 Denote by
 \begin{align}
 \eta : \widetilde{\mcT}_{\mcP_B/S} \migi  \mr{ad}(\mcP_G)
 \notag
 \end{align}
 the $\mcO_{C_S}$-linear  morphism 
  given by 
 $\partial \mapsto \partial - (\nabla_{\mcP_G} \circ \alpha_{\mcP_G})(\partial)$ for any local section $\partial \in \widetilde{\mcT}_{\mcP_B/S}$.
One verifies that $\eta$ is an isomorphism
and satisfies $\widetilde{\nabla}_{\mcP_B}^\mr{ad} =\nabla_{\mcP_G}^\mr{ad} \circ \eta$.
Hence,  the pair of morphisms $(\eta, \mr{id}_{\Omega_{C_S/S}\otimes \mr{ad}(\mcP_G)})$
   specifies
  an isomorphism $\eta^\bullet :  \mcK^\bullet [\widetilde{\nabla}^\mr{ad}_{\mcP_B}] \isom \mcK^\bullet [\nabla^\mr{ad}_{\mcP_G}]$.
This isomorphism  
 fits into the following isomorphism of sequences of complexes:
\begin{align}
\begin{CD}
0 @>>> \Omega_{C_S/S}^{\otimes 2}[-1]@>>> \mcK^\bullet [\widetilde{\nabla}^\mr{ad}_{\mcP_B}] @> {\alpha'}_{\mcP_B}^{\bullet}>>  \mcT_{C_S/S}[0]@>>> 0
\\
@. @V \wr V \mr{id} V @V \wr V \eta^\bullet V @V \wr V \mr{id} V @.
\\
0 @>>>\Omega_{C_S/S}^{\otimes 2}[-1] @>>> \mcK^\bullet [\nabla^\mr{ad}_{\mcP_G}] @>>> \mcT_{C_S/S} [0]@>>> 0,
\end{CD} \notag
\end{align}
where 
\begin{itemize}
\item
${\alpha'}_{\mcP_B}^{\bullet}$ denotes the morphism given by $\alpha'_{\mcP_B}$ (cf. (\ref{W51}));
\item
both the upper and lower second arrows arise from the composite
\begin{align}
  \Omega_{C_S/S}^{\otimes 2} \xrightarrow{\mr{id}\otimes \overline{\nabla}_{\mcP_G}^\sharp} \Omega_{C_S/S}\otimes \mr{ad}(\mcP_G)^2 \migiincl \Omega_{C_S/S}\otimes \mr{ad}(\mcP_G)^1; \notag
\end{align}
\item
the lower third arrow arises from the composite 
\begin{align}
 \mr{ad}(\mcP_G) \migisurj \mr{ad}(\mcP_G)/\mr{ad}(\mcP_G)^1 \xrightarrow{\overline{\nabla}^\flat_{\mcP_G}} \mcT_{C_S/S}.
\notag
\end{align}
\end{itemize}

By applying the functor $\mbR^1f_{S*}(-)$ to this diagram, we obtain
the following isomorphism of sequences of $\mcO_S$-modules: 
\begin{align} \label{WW1003} 
\begin{CD}
0 @>>> f_{S*}(\Omega_{C_S/S}^{\otimes 2})@> {\gamma'}^{\sharp}_{\mcP^\ind} >> \mbR^1 f_{S*}(\mcK^\bullet [\widetilde{\nabla}^\mr{ad}_{\mcP_B}]) @>{\gamma'}^{\flat}_{\mcP^\ind} >>  \mbR^1f_{S*}(\mcT_{C_S/S})@>>> 0
\\
@. @V \wr V \mr{id} V @V \wr V \mbR^1f_{S*}(\eta^\bullet) V @V \wr V \mr{id} V @.
\\
0 @>>>f_{S*}(\Omega_{C_S/S}^{\otimes 2})@>>\gamma^\sharp_{\mcP^\ind} > \mbR^1f_{S*}(\mcK^\bullet [\nabla^\mr{ad}_{\mcP_G}]) @>>\gamma^\flat_{\mcP^\ind} > \mbR^1 f_{S*}(\mcT_{C_S/S})@>>> 0.
\end{CD}
\end{align}
Since the lower horizontal sequence is   exact  (cf.  ~\cite{Mzk1}, Chap.\,I, \S\,2, Theorem 2.8 (1)), the upper sequence turns out to be exact.
Finally, we remark that the lower  sequence in (\ref{WW1003}) is obtained by taking cohomology of the differentials in the $E_1$-term of the spectral sequence
\begin{align} 
    E_{1}^{p,q} = \mbR^qf_{S*}(\mcK^p [{\nabla_{\mcP_G}^\mr{ad}}]) \Rightarrow \mbR^{p+q}f_{S*}(\mcK^\bullet [{\nabla_{\mcP_G}^\mr{ad}}]).   \notag
    \end{align}

\vspace{3mm}
\subsection{} \label{sss12}

Next, let us introduce the notion of an indigenous bundle {\it of canonical type}.
Write   $\pi^\dagger_{\mr{GL}_2} : \mcP_{\mr{GL}_2}^\dagger \migi C_S$
for the (right)  $\mr{GL}_2$-torsor over $C_S$ associated with the rank $2$ vector bundle ${^r}\mcD^{\leq 1}_{C_S/S}$.
 More precisely, $ \mcP_{\mr{GL}_2}^\dagger$ is the $C_S$-scheme representing the functor
\begin{align}
\mcI som_{\mcO_{C_S}} (\mcO_{C_S}^{\oplus 2}, {^r}\mcD^{\leq 1}_{C_S}) : \mfE \mft_{C_S} \migi \mfS \mfe \mft,\notag
\end{align}
classifying locally defined isomorphisms $\mcO_{C_S}^{\oplus 2} \isom  {^r}\mcD^{\leq 1}_{C_S}$,  where $\mfE \mft_{C_S}$ denotes the small \'{e}tale site on $C_S$.
Also, write
\begin{align}
\pi_G^\dagger : \mcP^\dagger_G \migi C_S\notag
\end{align}
for the $G$-torsor induced by $\mcP^\dagger_{\mr{GL}_2}$ via the change of structure group by the quotient $\mr{GL}_2 \migisurj G$.
Let us  consider $\mcO_{C_S}$ as an $\mcO_{C_S}$-submodule of $\mcO_{C_S}^{\oplus 2}$ via the injection $\mcO_{C_S} \migiincl \mcO_{C_S}^{\oplus 2}$ into the first factor.
 The subfunctor
of $\mcI som_{\mcO_{C_S}} (\mcO_{C_S}^{\oplus 2}, {^r}\mcD^{\leq 1}_{C_S})$
consisting of locally defined  isomorphisms $w : \mcO_{C_S}^{\oplus 2} \isom {^r}\mcD^{\leq 1}_{C_S}$ with $w (\mcO_{C_S}) \subseteq {^r}\mcD^{\leq 0}_{C_S/S}$ may be represented by a $B$-reduction
\begin{align}
\pi_B^\dagger : \mcP^\dagger_B \migi C_S\notag
\end{align}
of $\mcP_G^\dagger$.
The adjoint vector bundle $\mr{ad} (\mcP_G^\dagger)$ is canonically isomorphic to
the sheaf $\mcE nd_{\mcO_{C_S}}^0 ({^r}\mcD^{\leq 1}_{C_S/S})$ of  $\mcO_{C_S}$-linear  endomorphisms of ${^r}\mcD^{\leq 1}_{C_S/S}$ with vanishing trace.

Let  us take 
a scheme $U$ equipped with an \'{e}tale morphism $U \migi C_S$
 and a  section  $x \in \Gamma (U, \mcO_{C_S})$ such that $d x$ generates $\Omega_{U/S}$ ($= \Omega_{C_S/S} |_U$).
We shall refer to such a pair $(U, x)$ as a {\bf local chart} of $C_S$ relative to $S$.
  The   decomposition
${^r}\mcD^{\leq 1}_{C_S/S} |_U \isom  \mcO_{U} \oplus \mcO_{U} \cdot  \partial_x$ by means of $(U, x)$, where $\partial_x \in \Gamma (U, \mcT_{C_S/S})$ denotes the dual base  of $d x$, 
gives an isomorphism  
\begin{align}
\tau_{(U, x)} : \mcP_G^\dagger |_U \ (= U \times_{C_S} \mcP^\dagger_G)  \isom U \times_R G \notag
\end{align}
 of $G$-torsors which induces an isomorphism
 \begin{align} \label{eEEe2}
\widetilde{\tau}_{(U,x)}^\mr{ad} :  \widetilde{\mcT}_{\mcP^\dagger_G/S} |_U  \isom \mcT_{U/S} \oplus  (\mcO_{U} \otimes_R \mfg).
\end{align}

\vspace{3mm}
\bde  \label{dD01}
  We shall say that an indigenous bundle $\mcP^\circledast$ on $C_S/S$ is {\bf of canonical type}
 if it is of the form $\mcP^\circledast = (\mcP_B^\dagger, \nabla_{\mcP_G})$, where
 $\nabla_{\mcP_G}$ is an $S$-connection on $\mcP_G^\dagger$,  satisfying the following  condition:   
 for any local chart $(U, x)$ of $C_S$ relative to $S$, the restriction $\nabla_{\mcP_G} |_U$ of $\nabla_{\mcP_G}$ to $U$  may be  expressed, via (\ref{eEEe2}),  as the map $\mcT_{U/S} \migi \mcT_{U/S} \oplus (\mcO_U  \otimes_R \mfg)$ determined by 
 $1 \cdot \partial_x \mapsto \left(1 \cdot \partial_x, \begin{pmatrix} 0 & a \\ 1 & 0 \end{pmatrix}\right)$ for some $a \in \Gamma (U, \mcO_{C_S})$.
   \ede

\vspace{3mm}
\subsection{} \label{s1f3}

In \S\S\,\ref{s1f3}-\ref{sSf13},
 we  describe  indigenous bundles in terms of   differential operators between line bundles.
Recall that a {\it theta characteristic} (in other words,  a  {\it spin structure}) of $C_S/S$ is, by definition, a line bundle $\mcL$ on $C_S$ together with an isomorphism $\mcL^{\otimes 2} \isom \Omega_{C_S/S}$.
The curve $C_S/S$ necessarily admits, at least \'{e}tale locally on $S$, a theta characteristic
(cf. ~\cite{Wak}, Remark 2.2.1 (i)).

Let $\mcL$ be a theta characteristic of $C_S/S$.
We shall write
\begin{align} 
\mcD {\it iff}^{\leq 2}_\mcL & \,  := \mcH om_{\mcO_{C_S}} (\mcL^{\vee},  \mcL^\vee \otimes  \Omega_{C_S/S}^{\otimes 2}\otimes \mcD_{C_S/S}^{\leq  2}) 
 \  \left(\cong  \mcH om_{\mcO_{C_S}} (\mcL^{\vee},  \mcL^{\otimes 3}\otimes \mcD_{C_S/S}^{\leq  2}) \right),  \notag \\
  {^\circlearrowright \mcD} {\it iff}^{\leq 2}_\mcL
& \, := \mcH om_{\mcO_{C_S}} (\mcT_{C_S/S}^{\otimes 2} \otimes \mcL,  \mcD_{C_S/S}^{\leq 2}\otimes \mcL)
 \  \left(\cong \mcH om_{\mcO_{C_S}} (\mcL^{\otimes (-3)},  \mcD_{C_S/S}^{\leq 2}\otimes \mcL) \right), \notag \end{align}
  where $\mcL^\vee \otimes \Omega_{C_S/S}^{\otimes 2} \otimes \mcD_{C_S/S}^{\leq 2}$ ($\cong \mcL^{\otimes 3}\otimes \mcD_{C_S/S}^{\leq  2}$) and  $\mcD_{C_S/S}^{\leq 2}\otimes \mcL$ may be considered as being equipped with  the structures of $\mcO_{C_S}$-module defined   in \S\,\ref{sss13}.
A second order differential operator from $\mcL^\vee$ to $\mcL^\vee \otimes \Omega_{C_S/S}^{\otimes 2}$
is nothing but a global section of $\mcD {\it iff}^{\leq 2}_\mcL$.
By passing to the composite injection
\begin{align} 
\Omega_{C_S/S}^{\otimes 2} & \  \isom \mcH om_{\mcO_{C_S}} (\mcL^\vee, \mcL^\vee \otimes \Omega_{C_S/S}^{\otimes 2} \otimes  \mcO_{C_S})
  \left(= \mcH om_{\mcO_{C_S}} (\mcL^\vee, \mcL^\vee \otimes \Omega_{C_S/S}^{\otimes 2} \otimes \mcD^{\leq 0}_{C_S/S}) \right)
      \migiincl \mcD {\it iff}^{\leq 2}_\mcL, \notag
\end{align}
 we identify $\Gamma (C_S, \Omega_{C_S/S}^{\otimes 2})$  with a  submodule of $\Gamma (C_S, \mcD {\it iff}^{\leq 2}_\mcL)$.

Next, let us define 
\begin{align} 
\mcD {\it iff}^{\leq 2 \bigstar}_\mcL \notag
\end{align}
to be the subsheaf  of $\mcD {\it iff}^{\leq 2}_\mcL$ consisting of differential operators
such that the principal symbol is $1$ and the subprincipal symbol is $0$.
That is to say, 
a second order differential operator  $\mbD$ from $\mcL^\vee$ to $\mcL^\vee \otimes \Omega_{C_S/S}^{\otimes 2}$  lies in 
 $\mcD {\it iff}^{\leq 2 \bigstar}_\mcL$ if and only if
whenever we choose a local chart  $(U, x)$ of $C_S$  and a trivialization  $\mcL |_U \cong \mcO_{U} \cdot (dx)^{\otimes \frac{1}{2}}$,  the operator $\mbD$ may be   given, on $U$,  by assigning 
\begin{align}
(dx)^{\otimes (-\frac{1}{2})} \mapsto (dx)^{\otimes (-\frac{1}{2})} \otimes  (dx)^{\otimes 2} \otimes \partial^2_x + (dx)^{\otimes (-\frac{1}{2})} \otimes a \otimes 1\notag
\end{align}
 for some $a \in \Gamma (U, \Omega_{C_S/S}^{\otimes 2})$.
For each $\mbD \in \Gamma (C_S, \mcD {\it iff}^{\leq 2 \bigstar}_\mcL)$ and $A \in \Gamma (C_S, \mcD {\it iff}^{\leq 2}_\mcL)$,
 the sum $\mbD + A$ lies in $\Gamma (C_S, \mcD {\it iff}^{\leq 2 \bigstar}_\mcL)$
 if and only if $A \in \Gamma (C_S, \Omega_{C_S/S}^{\otimes 2})$ ($\subseteq \Gamma (C_S, \mcD {\it iff}^{\leq 2}_\mcL)$).
 Thus, {\it $\Gamma (C_S, \mcD {\it iff}^{\leq 2 \bigstar}_\mcL)$  admits a canonical structure of   $\Gamma (C_S, \Omega_{C_S/S}^{\otimes 2})$-torsor}.

\vspace{3mm}
\subsection{} \label{sSf13}

Now, let   $\mbD$ be an element of $\Gamma (C_S, \mcD {\it iff}^{\leq 2 \bigstar}_\mcL)$.
Denote by ${^\circlearrowright \mbD} : \mcT_{C_S/S}^{\otimes 2} \otimes \mcL \migi \mcD^{\leq 2}_{C_S/S} \otimes \mcL$ the morphism defined  as the image of $\mbD$  via the composite
\begin{equation} 
\label{GL12}
{\mcD} {\it iff}^{\leq 2}_\mcL  \isom   \mcL^{\otimes 3} \otimes \mcD_{C_S/S}^{\leq 2} \otimes \mcL \isom {^\circlearrowright \mcD} {\it iff}^{\leq 2}_\mcL
\end{equation}
of two natural isomorphisms arising from the definitions of ${\mcD} {\it iff}^{\leq 2}_\mcL$ and ${^\circlearrowright \mcD} {\it iff}^{\leq 2}_\mcL$.
Since  $\mbD \in \Gamma (C_S, \mcD {\it iff}^{\leq 2 \bigstar}_\mcL)$,  the composite 
\begin{align} \label{eEEe15}
\mcD^{\leq 1}_{C_S/S} \otimes \mcL \migiincl \mcD_{C_S/S} \otimes \mcL \migisurj (\mcD_{C_S/S} \otimes \mcL) / \langle \mr{Im} ({^\circlearrowright \mbD})\rangle
\end{align}
is an isomorphism,
 where  $(\mcD_{C_S/S} \otimes \mcL)/ \langle \mr{Im} ({^\circlearrowright \mbD}) \rangle$ denotes the quotient  of $\mcD_{C_S/S} \otimes \mcL$ by the left $\mcD_{C_S/S}$-submodule generated by $\mr{Im} ({^\circlearrowright \mbD})$.
The structure of left $\mcD_{C_S/S}$-module on $\mcD_{C_S/S}^{\leq 1} \otimes \mcL$ transposed  from $(\mcD_{C_S/S} \otimes \mcL)/ \langle \mr{Im} ({^\circlearrowright \mbD}) \rangle$ via (\ref{eEEe15})  corresponds (cf. \S\,\ref{sss13}) to
an $S$-connection 
\begin{align}
\nabla_\mbD : \mcD^{\leq 1}_{C_S/S} \otimes  \mcL \migi \Omega_{C_S/S} \otimes (\mcD^{\leq 1}_{C_S/S} \otimes \mcL)\notag
\end{align}
 on $\mcD^{\leq 1}_{C_S/S} \otimes \mcL$.
 Note that (since the subprincipal symbol of $\mbD$ is $0$) the $S$-connection $\mr{det} (\nabla_\mbD)$ on the determinant $\mr{det} (\mcD^{\leq 1}_{C_S/S} \otimes \mcL)$
  induced by $\nabla_\mbD$  coincides,  
  via 
 $\mr{det} (\mcD^{\leq 1}_{C_S/S} \otimes \mcL) \cong \mcT_{C_S/S} \otimes \mcL^{\otimes 2} \cong \mcO_{C_S}$, with the universal derivation $d : \mcO_{C_S} \migi \Omega_{C_S/S}$.
 The $G$-torsor associated 
  with $\mcD^{\leq 1}_{C_S/S} \otimes \mcL$ via projectivization  is canonically isomorphic to $\mcP^\dagger_G$.
 Hence, $\nabla_\mbD$ yields   an $S$-connection $\nabla_{\mbD, G}$ on $\mcP^\dagger_G$.
 It follows from the various definitions involved that the pair
 \begin{align} 
 \mcP^{\circledast \dagger}_\mbD := (\mcP^\dagger_B, \nabla_{\mbD, G})\notag
 \end{align} 
 forms an indigenous bundle on $C_S/S$ {\it of canonical type}.


\vspace{3mm}
\bpr \label{pP01}
\begin{itemize}
\item[(i)]
Suppose that there exists   a theta characteristic $\mcL$ of $C_S/S$.
Then, 
the assignment $\mbD \mapsto \mcP^{\circledast \dagger}_{\mbD}$ constructed above  defines 
a bijective correspondence between the set $\Gamma (C_S, \mcD {\it iff}^{\leq 2 \bigstar}_\mcL)$ and the set of isomorphism classes of indigenous bundles on $C_S/S$.
\item[(ii)]
For any indigenous bundle $\mcP^\circledast$ on $C_S/S$, 
there exists a unique  indigenous bundle  $\mcP^{\circledast \dagger}$ on $C_S/S$ of canonical type
which is isomorphic to $\mcP^\circledast$.
\end{itemize}
 \epr
\begin{proof}
Let us consider   assertion (i).
First, we shall prove  the injectivity of the assignment $\mbD \mapsto \mcP^{\circledast \dagger}_\mbD$.
Let $\mbD_1$ and $\mbD_2$ be 
 elements of
 $\Gamma (C_S, \mcD {\it iff}^{\leq 2 \bigstar}_\mcL)$,  and suppose    that 
 there exists an isomorphism $\xi : \mcP^{\circledast \dagger}_{\mbD_1} \isom \mcP^{\circledast \dagger}_{\mbD_2}$.
Let  $(U, x)$ be  a local chart of $C_S$ relative to $S$, which gives an identification $\tau_{(U, x)}  : \mcP^\dagger_{G} |_U \isom U \times_R G$.
  After possibly replacing $U$ with its \'{e}tale covering, one may find   an element $R := \begin{pmatrix} b & c \\ 0 & \frac{1}{b} \end{pmatrix} \in \mr{SL}_2 (U)$  
 such that the restriction of $\xi$ to $U$ may be described, via $\tau_{(U, x)}$,  as  the assignment   $v \mapsto \overline{R} \cdot v$ (for any  $v \in B (U)$), where  $\overline{R}$ denotes   the image of $R$ via the quotient $\mr{SL}_2 \migisurj G$.
Here, observe that the restrictions $\nabla_{\mbD_1} |_U$, $\nabla_{\mbD_2} |_U$ of the $S$-connections $\nabla_{\mbD_1}$, $\nabla_{\mbD_2}$ may be expressed, via $\widetilde{\tau}_{(U, x)}^\mr{ad}$ (cf. (\ref{eEEe2})),  as the maps  determined, repectively, by 
\begin{align}
1 \cdot \partial_x \mapsto \left(1 \cdot \partial_x, \begin{pmatrix} 0 & a_1 \\ 1 & 0 \end{pmatrix}\right), 
 \ \ \ 1 \cdot \partial_x \mapsto \left(1 \cdot \partial_x, \begin{pmatrix} 0 & a_2 \\ 1 & 0 \end{pmatrix}\right) \notag
 \end{align}
  for some $a_1, a_2  \in \Gamma (U, \mcO_{C_S})$.
Since $\xi$ is compatible with the respective $S$-connections $\nabla_{\mbD_1}$ and $\nabla_{\mbD_2}$, the following equalities hold in  $\Gamma (U,   \mcO_{C_S} \otimes_R \mfg \mfl_2)$:
\begin{align} 
\begin{pmatrix} 0 & a_2 \\ 1 & 0 \end{pmatrix}  & =
 R  \cdot \begin{pmatrix} 0 & a_1 \\ 1 & 0 \end{pmatrix}  \cdot R^{-1}
      +  
      R \cdot d (R^{-1})  \ \left(=   \begin{pmatrix} \frac{-b' + c}{b} & a_1b^2 -c^2 -bc' + b' c \\ \frac{1}{b^2} & \frac{b'-c}{b} \end{pmatrix}\right). \notag
\end{align}
It follows   that
($b^2 =1$, $c =0$, and) $a_1 = a_2 $.
In particular, by applying the above discussion to various local charts $(U, x)$,  we obtain the equality 
$\mbD_1 = \mbD_2$.
This completes the proof of the  injectivity.

Next, let us prove the surjectivity of the assignment $\mbD \mapsto \mcP^{\circledast \dagger}_\mbD$.
Let $\mcP^\circledast := (\mcP_B, \nabla_{\mcP_G})$ be an indigenous bundle on $C_S/S$.
It follows from Proposition  \ref{R002} that, by means of  descent with respect to \'{e}tale morphisms,  
we are    always  free to replace $S$ by any  \'{e}tale covering of $S$.
Hence, 
we may assume that there exists
 a rank $2$ vector bundle $\mcV$ on $C_S$ {\it with trivial determinant} 
such that the $G$-torsor associated with  $\mcV$ via 
  projectivization is isomorphic to $\mcP_G$.
The $B$-reduction $\mcP_B$ of $\mcP_G$ determines a line subbundle $\mcN$ of $\mcV$, which satisfies that $\mcV/\mcN \cong \mcN^\vee$.
Here, 
write  
$\mcP_{\mr{GL}_2}$ for  the $\mr{GL}_2$-torsor corresponding to
$\mcV$ and 
 $\mcP_{\mbG_m}$ for 
 the $\mbG_m$-torsor
induced from  $\mcP_{\mr{GL}_2}$ via the change of structure group by 
the determinant map $\mr{det} : \mr{GL}_2 \migi \mbG_m$.
Since   the map $(q, \mr{det}) : \mr{GL}_2 \isom G \times_R \mbG_m$, where $q$ denotes the natural quotient $\mr{GL}_2 \migisurj G$,   is \'{e}tale, 
the natural morphism
 \begin{align}
 \widetilde{\mcT}_{\mcP_{\mr{GL}_2}/S} \isom \widetilde{\mcT}_{\mcP_G/S} \times_{\alpha_{\mcP_G}, \mcT_{C_S/S}, \alpha_{\mcP_{\mbG_m}}} \widetilde{\mcT}_{\mcP_{\mbG_m}/S} \notag
 \end{align}
is an isomorphism.
By this isomorphism,  the pair $(\nabla_{\mcP_G}, d)$ of    the $S$-connection $\nabla_{\mcP_G}$ on $\mcP_G$ and the trivial $S$-connection $d$ on $\mcO_{C_S}$
determines an $S$-connection $\nabla_{\mcV}$ on $\mcV$. 
For each $j = 0, 1, 2$, we shall  write $\zeta^j$ for  the composite
\begin{align} 
\zeta^j : \mcD^{\leq j}_{C_S/S} \otimes \mcN \migiincl \mcD_{C_S/S} \otimes \mcV \xrightarrow{\widehat{\nabla}_{\mcV}} \mcV,\notag
\end{align}
where the first arrow arises from the inclusions  $\mcD^{\leq j}_{C_S/S} \migiincl \mcD_{C_S/S}$ and $\mcN \migiincl \mcV$.
In particular, $\zeta^0$ coincides  with
the inclusion $\mcN \migi \mcV$ under  the natural identification $\mcD^{\leq 0}_{C_S/S} \otimes \mcN = \mcN$.
 It follows from the definition of an indigenous bundle that $\zeta^1$ is an isomorphism and 
$\zeta^2$ is surjective.
Moreover,  $\zeta^1$ induces, by taking determinants,  an isomorphism 
 \begin{align} \label{eEEe19}
 \mcT_{C_S/S} \otimes \mcN^{\otimes 2} \ \left(\cong \mr{det} (\mcD^{\leq 1}_{C_S/S} \otimes \mcN)\right) \isom \left(\mr{det} (\mcV) \cong \right) \ \mcO_{C_S}.
 \end{align}
That is to say, the line bundle $\mcN$ together with an isomorphism $\mcN^{\otimes 2} \isom \Omega_{C_S/S}$ induced by (\ref{eEEe19}) specifies  a theta characteristic of  $C_S/S$.
The square of $\mcL' := \mcL \otimes \mcN^\vee$ is trivial, so
there exists  a unique $S$-connection  $\nabla_{\mcL'}$  on $\mcL'$ whose square coincides with $d$.
By tensoring $(\mcV, \nabla_\mcV)$ with $(\mcL', \nabla_{\mcL'})$,
we may assume that the vector bundle $\mcV$ was taken to satisfy $\mcN = \mcL$.

Let us consider the short exact sequence
\begin{align} \label{eEEe5}
0 \migi \mcD^{\leq 1}_{C_S/S} \otimes \mcL \migi \mcD^{\leq 2}_{C_S/S} \otimes \mcL \migi \mcT^{\otimes 2}_{C_S/S} \otimes \mcL \migi 0
\end{align}
(cf. (\ref{GL7})).
The composite $(\zeta^1)^{-1} \circ \zeta^2 : \mcD^{\leq 2}_{C_S/S} \otimes \mcL \migisurj \mcD^{\leq 1}_{C_S/S} \otimes \mcL$ specifies    a split surjection of  (\ref{eEEe5}).
This split surjection determines a split injection $\mcT^{\otimes 2}_{C_S/S} \otimes \mcL \migi \mcD^{\leq 2}_{C_S/S} \otimes \mcL$ of (\ref{eEEe5}),
and hence, determines a differential operator $\mbD : \mcL^\vee\migi \mcL^\vee \otimes \Omega_{C_S/S}^{\otimes 2} \otimes \mcD^{\leq 2}_{C_S/S}$ via (\ref{GL12}).
The isomorphism $\zeta^1$ is verified to be compatible with the respective $S$-connections $\nabla_\mbD$ and $\nabla_\mcV$, as well as with the respective filtrations $\mcD^{\leq 1}_{C_S/S} \otimes \mcL  \supseteq \mcD^{\leq 0}_{C_S/S} \otimes \mcL \supseteq 0$ and $\mcV \supseteq \mcL \supseteq 0$.
Consequently, $\zeta^1$ determines an isomorphism $\mcP_\mbD^{\circledast\dagger} \isom \mcP^\circledast$ of indigenous bundles.
This implies that the assignment $\mbD \mapsto \mcP_\mbD^{\circledast\dagger}$ is surjective, and we finish  the proof of  assertion (i).

Next, let us consider assertion (ii).
Since indigenous bundles (of canonical type) may be constructed by means of descent with respect to \'{e}tale morphisms, 
we are  
always free to replace $S$ with its \'{e}tale covering (cf.  Proposition \ref{R002}).
In particular, the problem is  reduced to the case where 
  $C_S/S$ admits   a theta characteristic.
  Hence,  assertion (ii)  follows from assertion (i).
\end{proof}

\vspace{3mm}
\subsection{} \label{s13}

Let us introduce notations concerning moduli functors classifying indigenous bundles.
Denote by 
\begin{align} 
 \mfS_{g, R} : \mfS \mfc \mfh_{/\mfM_{g, R}} \migi \mfS \mfe \mft  \notag
\end{align}
the $\mfS \mfe \mft$-valued functor on $\mfS \mfc \mfh_{/\mfM_{g,R}}$  which, to any   
object $T \migi \mfM_{g,R}$ of $\mfS \mfc \mfh_{/\mfM_{g,R}}$,
 assigns the set of isomorphism classes of indigenous bundles on the curve $f_T : C_T \migi T$.
Also, denote by
\begin{align}
\pi_{g,R}^\mfS : \mfS_{g, R} \migi \mfM_{g, R} \notag
\end{align}
the natural projection.
Given an object $S \migi \mfM_{g, R}$ of $\mfS \mfc \mfh_{/ \mfM_{g,R}}$, we obtain 
\begin{align} 
  \mfS_S := \mfS_{g, R} \times_{\pi_{g,R}^\mfS, \mfM_{g,R}} S, \notag
  \end{align}
which has the projection
\begin{align} 
   \pi^\mfS_S : \mfS_S \migi S. \notag
   \end{align}
 
In what follows, let us consider 
   a natural affine structure on $\mfS_S$
by means of modular interpretation.
Let $S$ be as above.
Also, 
let $\mcP^\ind := (\mcP_B, \nabla_{\mcP_G})$ be an indigenous bundle on  $C_S/S$ 
and
$A$ an element of $\Gamma (C_S, \Omega_{C_S/S}^{\otimes 2})$, or equivalently a global  section of the projection $\mbA(f_{S*}(\Omega_{C_S/S}^{\otimes 2})) \migi S$.
By Propositions \ref{R002} and \ref{pP01} (ii),
there exist a unique indigenous bundle $\mcP^{\circledast \dagger} := (\mcP^\dagger_B, \nabla^\dagger_{\mcP_G})$ on $C_S/S$ of canonical type and a unique isomorphism $\mcP^\circledast \isom \mcP^{\circledast \dagger}$.
By passing to the composite 
\begin{align}  
\Omega_{C_S/S}^{\otimes 2}  & \xrightarrow{\mr{id}\otimes \overline{\nabla}^{\dagger \sharp}}  \Omega_{C_S/S}\otimes \mr{ad}(\mcP_G^\dagger)^2 
 \stackrel{}{\migiincl}  \Omega_{C_S/S} \otimes  \widetilde{\mcT}_{\mcP^\dagger_G/S}
    \isom    \mcH om_{\mcO_{C_S}} (\mcT_{C_S/S},  \widetilde{\mcT}_{\mcP_G^\dagger/S}), \notag 
\end{align}
one may identify $A$ with  an $\mcO_{C_S}$-linear morphism $\mcT_{C_S/S} \migi \widetilde{\mcT}_{\mcP^\dagger_\G/S}$.
Hence, the sum $\nabla^\dagger_{\mcP_G} +A : \mcT_{C_S/S} \migi \widetilde{\mcT}_{\mcP_G^\dagger/S}$ makes sense and 
 specifies  an $S$-connection on $\mcP_G^\dagger$.
Moreover, it follows from the definition of an indigenous bundle (of canonical type) that
the pair 
\begin{align} 
\mcP^{\ind \dagger}_{+A}:=(\mcP_B^\dagger, \nabla^\dagger_{\mcP_G}+A) \notag
\end{align}
 forms an indigenous bundle on $C_S/S$ of canonical type.
Notice that  if  there exist a theta characteristic $\mcL$ of $C_S/S$ and  an element $\mbD$ of 
$\Gamma (C_S, \mcD {\it iff}_\mcL^{\leq 2 \bigstar})$ with 
 $\mcP^{\circledast \dagger} = \mcP^{\circledast \dagger}_\mbD$ (cf. Proposition \ref{pP01} (ii)), 
 then  
 $\mcP^{\circledast \dagger}_{+A}$ may be described as   
  $\mcP^{\circledast \dagger}_{+A} = \mcP^{\circledast \dagger}_{\mbD + A}$.
The assignment $(\mcP^\ind, A) \mapsto \mcP^{\ind \dagger}_{+A}$
is well-defined and  functorial  with respect to $S$.
Therefore, it
 determines an action 
\begin{align} 
   \mfS_S \times_S \mbA(f_{S*}(\Omega_{C_S/S}^{\otimes 2})) \migi \mfS_S \notag
   \end{align}
of the relative affine space $\mbA(f_{S*}(\Omega_{C_S/S}^{\otimes 2}))$ on $\mfS_S$.
By
 Proposition \ref{pP01} (i) and  (ii), the following proposition holds.

\vspace{3mm}
\bpr [cf.  ~\cite{Mzk1}, Chap.\,I, \S\,2, Corollary 2.9] \label{p01}
The functor $\mfS_S$ may be represented by an $\mbA (f_{S*}(\Omega_{C_S/S}^{\otimes 2}))$-torsor over $S$ with respect to the  $\mbA (f_{S*}(\Omega_{C_S/S}^{\otimes 2}))$-action just discussed.
In particular, if, moreover, $S$ is a geometrically connected  smooth Deligne-Mumford stack over $R$ of relative dimension $n$,
then the functor $\mfS_S$ may be represented by a geometrically connected smooth Deligne-Mumford stack  over $R$ of relative dimension $n + 3g-3$.
 \epr

\vspace{3mm}
\subsection{} \label{Ws35}

In this subsection, we shall  consider a  cohomological expression of 
 the deformation space of an indigenous bundle.
Let $S$  and $\mcP^\ind$ be as above.

Until just before Proposition \ref{p02}, we impose the  assumption that  $S$ is affine.
Then, the sequence (\ref{Ww200}) reads 
\begin{align} \label{Ww600}
0 \longmigi \mbH^1 (\mcK^\bullet[\nabla^\mr{ad}_{\mcP_G}])
\longmigi \mbH^1 (\mcK^\bullet[\widetilde{\nabla}^\mr{ad}_{\mcP_G}])
\xrightarrow{\alpha_{\mcP_G}^\mbH} H^1 (C_S, \mcT_{C_S/S}) \longmigi 0,
\end{align}
where, given a complex  $\mcK^\bullet$, we shall denote by 
$\mbH^1 (\mcK^\bullet)$ its $1$-st hypercohomology group.
In what follows,  we shall  write $R_\epsilon := R [\epsilon]/(\epsilon^2)$, and
  denote the base-changes  to $R_\epsilon$ of objects over $R$  by means of a subscripted $\epsilon$.
Write $v :  S  \migi \mfS_{g, R}$ for 
the 
classifying morphism of $(C_S, \mcP^\circledast)$
and $\overline{v} : S \migi \mfM_{g, R}$
for the classifying  morphism of  $C_S$, i.e., $\overline{v} := \pi^{\mfS}_{g, R} \circ v$.

The tangent space
\begin{align} 
T_{C_S} := \Gamma (S, \overline{v}^*(\mcT_{\mfM_{g, R}/R})) \notag
\end{align}
of $\mfM_{g, R}/R$ at $\overline{v}$
may be identified with 
the 
deformation space of
  the curve $C_S$ over $S_\epsilon$.
We shall denote by
\begin{align} 
T_{C_S, \mcP_G, \nabla_{\mcP_G}} \notag
\end{align}
the 
 deformation space  of $(C_S, \mcP_G, \nabla_{\mcP_G})$ (as a data consisting of  a curve  and a flat $G$-torsor over it) over $S_\epsilon$;
it has  the subspace 
\begin{align} 
T_{\mcP_G, \nabla_{\mcP_G}} \notag
\end{align}
 classifying 
 deformations 
 whose underlying curves are the trivial deformation.
By  arguments similar to the arguments in ~\cite{Chen}, \S\,4 (e.g., the proof of Proposition 4.4),
there are canonical isomorphisms
\begin{align} 
& t_{C_S} \hspace{12mm}: T_{C_S} \isom H^1 (C_S, \mcT_{C_S/S}), \notag \\
& t_{C_S, \mcP_G, \nabla_{\mcP_G}} :  T_{C_S, \mcP_G, \nabla_{\mcP_G}}   \isom
\mbH^1 (\mcK^\bullet [\widetilde{\nabla}^\mr{ad}_{\mcP_G}]),  \notag \\
& t_{\mcP_G, \nabla_{\mcP_G}} \hspace{4.5mm}:  T_{\mcP_G, \nabla_{\mcP_G}}  \isom \mbH^1 (\mcK^\bullet [\nabla^\mr{ad}_{\mcP_G}]) \notag
\end{align}
(cf. \S\,\ref{Ws3005} for their precise constructions)  making the following diagram commute:
\begin{align}\label{Egw110}
\begin{CD}
T_{\mcP_G, \nabla_{\mcP_G}}
@> \mr{incl.} >> T_{C_S, \mcP_G, \nabla_{\mcP_G}}
 @> 
 > >
T_{C_S}\\
@V \wr V t_{\mcP_G, \nabla_{\mcP_G}} V @V \wr V t_{C_S, \mcP_G, \nabla_{\mcP_G}} V @V \wr V t_{C_S} V
\\
\mbH^1 (\mcK^\bullet [\nabla^\mr{ad}_{\mcP_G}]) @>>> \mbH^1 (\mcK^\bullet [\widetilde{\nabla}^\mr{ad}_{\mcP_G}]) 
 @> > \alpha^\mbH_{\mcP_G}
  > H^1 (C_R, \mcT_{C_S/S}),
\end{CD} 
\end{align}
where the right-hand  arrow in the upper horizontal sequence is obtained  by forgetting the data of  deformations of $(\mcP_G, \nabla_{\mcP_G})$ and the lower horizontal arrows are the morphisms in (\ref{Ww600}).
Indeed, $t_{C_S}$ is nothing but
the Kodaira-Spencer map of $C_S/S$
(cf. ~\cite{Sch}, \S\,10, p.122).

Next, denote by 
\begin{align}
T_{C_S, \mcP^\circledast} \notag
\end{align}
the deformation space of $(C_S, \mcP_B, \nabla_{\mcP_G})$ (as a data consisting of a curve, a $B$-torsor over it, and an $S$-connection on the $G$-torsor associated to this $B$-torsor) over $S_\epsilon$.
Just as in the case of $T_{C_S, \mcP_G, \nabla_{\mcP_G}}$, 
there exists a canonical bijection
\begin{align} 
t_{C_S, \mcP^\circledast} : T_{C_S, \mcP^\circledast} \isom \mbH^1 (\mcK^\bullet [\widetilde{\nabla}^\mr{ad}_{\mcP_B}]), 
\notag
\end{align}
which
makes the following square diagram commute:
\begin{align} \label{Ww306}
\begin{CD}
T_{C_S, \mcP^\circledast}
@>>> T_{C_S, \mcP_G, \nabla_{\mcP_G}}
\\
@Vt_{C_S, \mcP^\circledast} V \wr V @V \wr V t_{C_S, \mcP_G, \nabla_{\mcP_G}} V
\\
\mbH^1 (\mcK^\bullet [\widetilde{\nabla}^\mr{ad}_{\mcP_B}]) 
@>> \mr{incl.}>
\mbH^1 (\mcK^\bullet [\widetilde{\nabla}^\mr{ad}_{\mcP_G}]),
\end{CD}
\end{align}
where the upper horizontal arrow arises from the  change of structure group by $B \migiincl G$.
Moreover,  
$T_{C_S, \mcP^\circledast}$ may be identified with  
the tangent space $\Gamma (S, v^*(\mcT_{\mfS_{g,R}/R}))$ of $\mfS_{g,R}/R$ at $v$, i.e.,
the deformation space of $(C_S, \mcP^\circledast)$ (as a pair of a curve and an indigenous bundle on it) over $S_\epsilon$.
Indeed, let
  $(C_{S}^v, \mcP_{B}^v, \nabla_{\mcP_{G}^v})$  be the deformation 
  classified by an element $v$ of 
  $T_{C_S, \mcP^\circledast}$, where we shall  write $\mcP_G^v := \mcP_B^v \times^B G$ and $\widetilde{\iota}^v : \widetilde{\mcT}_{\mcP_B^v/S_\epsilon} \migiincl \widetilde{\mcT}_{\mcP_G^v/S_\epsilon}$.
Then, the  composite
 \begin{align}
 \overline{\nabla}_{\mcP_{G}^v} :
 \mcT_{C_{S}^v/S_\epsilon}
  \xrightarrow{\nabla_{\mcP_{G}^v}} 
   \widetilde{\mcT}_{\mcP_{G}^v/S_\epsilon}
    \migisurj 
    \widetilde{\mcT}_{\mcP_{G}^v/S_\epsilon} / \widetilde{\iota}^v(\widetilde{\mcT}_{\mcP_{B}^v/S_\epsilon})\notag
 \end{align}
 is an isomorphism 
 since it becomes the isomorphism $\overline{\nabla}_{\mcP_G}$ when     restricted  to the  closed subscheme  $C_S \ \left(\subseteq C_{S}^v \right)$.  
 It follows that the pair $(\mcP_{B}^v, \nabla_{\mcP_G^v})$ forms an indigenous bundle on $C_{S}^v/S_\epsilon$, as desired.

On the other hand, 
the structure of  $\mbA (f_{S*}(\Omega_{C_S/S}^{\otimes 2}))$-torsor on $\mfS_{S}$ (cf. 
 Proposition \ref{p01}) yields a canonical bijection
\begin{align} 
  t_{\nabla_{\mcP_G}} : 
  T_{\nabla_{\mcP_G}} \ \left(:= \Gamma (S, v^*(\mcT_{\mfS_{g, R}/\mfM_{g,R}})) \right)
  \isom  \Gamma (C_S, \Omega_{C_S/S}^{\otimes 2}). \notag\end{align}
We have the following morphism of short exact sequences:
\begin{align}
\begin{CD}
0 @>>> T_{\nabla_{\mcP_G}}  @>>> T_{C_S, \mcP^\circledast} @>>> T_{C_S}@>>> 0
\\
@. @V \wr V t_{\nabla_{\mcP_G}}V @V \wr V t_{C_S, \mcP^\circledast} V @V \wr Vt_{C_S}V @.
\\
0 @>>> \Gamma (C_S, \Omega_{C_S/S}^{\otimes 2})  @>> {\gamma'}^\sharp_{\mcP^\circledast} > \mbH^1 (\mcK^\bullet [\widetilde{\nabla}^\mr{ad}_{\mcP_B}])@>> {\gamma'}^\flat_{\mcP^\circledast} > H^1 (C_S, \mcT_{C_S/S})@>>> 0,
\end{CD} \notag 
\end{align}
where 
the upper horizontal sequence is obtained by differentiating the smooth morphism  $\pi^\mfS_{g,R} : \mfS_{g, R} \migi \mfM_{g, R}$.

Finally, we remark that the various  isomorphisms obtained above are  functorial, in the natural sense, with respect to $S$.
Hence, by combining  (\ref{WW1003}) and the above isomorphism of short exact sequences for various affine schemes  $S$,
we obtain the following proposition.

\vspace{3mm}
\bpr  \label{p02}
Let us keep the above notation (but  the affineness assumption on $S$ are not imposed now).
Then, there exists a canonical isomorphism of short exact sequences of $\mcO_S$-modules:
\begin{align} \label{WW1004}
\begin{CD}
0 @>>> v^*(\mcT_{\mfS_{g, R}/\mfM_{g, R}})  @>>> v^*(\mcT_{\mfS_{g,R}/R}) @>>> \overline{v}^*(\mcT_{\mfM_{g, R}/R})@>>> 0
\\
@. @V \wr VV @V \wr V V @V \wr VV @.
\\
0 @>>> f_{S*}(\Omega_{C_S/S}^{\otimes 2})  @>> {\gamma}^\sharp_{\mcP^\circledast} >
 \mbR^1 f_{S*} (\mcK^\bullet [\nabla^\mr{ad}_{\mcP_G}])@>> {\gamma}^\flat_{\mcP^\circledast} > \mbR^1 f_{S*}(S, \mcT_{C_S/S})@>>> 0,
\end{CD} 
\end{align}
where 
the upper horizontal sequence is obtained by differentiating the smooth morphism  $\pi^\mfS_{g,R} : \mfS_{g, R} \migi \mfM_{g, R}$.
\epr

\vspace{3mm}
\subsection{} \label{Ws3005}

 We shall describe explicitly the deformations of data discussed above  in terms of \v{C}ech cohomology.
 To do this,  it suffices, by taking account of the commutative diagrams (\ref{Egw110}) and  (\ref{Ww306}),  to consider  
   (the inverse of) the bijective correspondence $t_{C_S, \mcP_G, \nabla_{\mcP_G}}$. 

Let us take an affine open covering $\mfU := \{ U_\alpha \}_{\alpha \in I}$ of $C_S$,   where $I$ is an index set.
We shall write $I_2$ for the set of pairs $(\alpha, \beta) \in I \times I$ with $U_{\alpha \beta} := U_\alpha \cap U_\beta \neq \emptyset$.
One may calculate $\mbH^1 (\mcK^\bullet[\widetilde{\nabla}^\mr{ad}_{\mcP_G}])$ as the total cohomology   of the \v{C}ech double complex   $\mr{Tot}^\bullet (\check{C}^\bullet (\mfU, \mcK^\bullet [\widetilde{\nabla}^\mr{ad}_{\mcP_G}]))$
 associated to $\mcK^\bullet [\widetilde{\nabla}^\mr{ad}_{\mcP_G}]$.
Each element $v$ of $\mbH^1 (\mcK^\bullet[\widetilde{\nabla}^\mr{ad}_{\mcP_G}])$  may be given by
a $1$-cocycle of $\mr{Tot}^\bullet (\check{C}^\bullet (\mfU, \mcK^\bullet [\widetilde{\nabla}^\mr{ad}_{\mcP_G}]))$, i.e.,
a collection of data
\begin{align}\label{W105}
 v = (\{ a_{\alpha \beta} \}_{\alpha, \beta}, \{ b_\alpha \}_{\alpha})
\end{align}
consisting of a \v{C}ech $1$-cocycle 
$\{ a_{\alpha\beta} \}_{\alpha, \beta} \in \check{C}^1 (\mfU, \widetilde{\mcT}_{\mcP_G/S})$, where $a_{\alpha\beta} \in \Gamma (U_{\alpha \beta}, \widetilde{\mcT}_{\mcP_G/S})$,  and a \v{C}ech $0$-cochain $\{b_\alpha \}_\alpha \in \check{C}^0 (\mfU, \Omega_{C_S/S} \otimes \mr{ad}(\mcP_G))$, where $b_\alpha \in \Gamma (U_\alpha, \Omega_{C_S/S} \otimes \mr{ad}(\mcP_G))$ $= \mr{Hom}_{\mcO_{U_\alpha}}(\mcT_{U_\alpha/S}, \mr{ad}(\mcP_G)|_{U_\alpha})$, 
which agree under $\widetilde{\nabla}^\mr{ad}_{\mcP_G}$ and the \v{C}ech coboundary map.
The elements in $\mbH^1 (\mcK^\bullet [\nabla^\mr{ad}_{\mcP_G}])$ (resp., $\mbH^1 (\mcK^\bullet [\widetilde{\nabla}^\mr{ad}_{\mcP_B}])$) may be represented by $v$ as above such that
$\{ a_{\alpha \beta} \}_{\alpha, \beta} \in \check{C}^1 (\mfU, \mr{ad}(\mcP_G))$ (resp., $\{ a_{\alpha \beta} \}_{\alpha, \beta} \in \check{C}^1 (\mfU, \widetilde{\mcT}_{\mcP_B/S})$).
The $S_\epsilon$-schemes $U_{\alpha, \epsilon} \ \left(:= U_{\alpha} \times_S S_\epsilon \right)$ for various $\alpha \in I$ may be glued together by means of the isomorphisms
\begin{align} 
\tau^v_{C_S, \alpha\beta} : =\mr{id}_{U_{\alpha \beta, \epsilon}} +\epsilon \cdot  \alpha_{\mcP_G} (a_{\alpha \beta}) : 
U_{\beta, \epsilon} |_{U_{\alpha \beta}} \isom  U_{\alpha, \epsilon} |_{U_{\alpha \beta}} \notag
\end{align}
 for $(\alpha, \beta) \in I_2$.
 The resulting $S_\epsilon$-scheme, which we denote by  $C_{S}^v$,
 specifies  the deformation corresponding to $\alpha_{\mcP_G}^\mbH (v) \ \left(\in H^1(C_S, \mcT_{C_S/S})\right)$  via
 $t_{C_S}$.
 Moreover,  the flat $G$-torsors  
 $(\mcP_{G, \epsilon} |_{U_\alpha}, \nabla_{\mcP_G, \epsilon} |_{U_\alpha}+ \epsilon \cdot b_\alpha)$
 may be glued together by means of the isomorphisms
 \begin{align}
\tau^v_{\mcP_G, \alpha \beta} := \mr{id}_{\mcP_{G, \epsilon} |_{U_{\alpha\beta}}} + \epsilon \cdot a_{\alpha \beta}
 :
   (\mcP_{G, \epsilon} |_{U_\beta}, \nabla_{\mcP_G, \epsilon} |_{U_\beta}+ \epsilon  \cdot b_\beta) |_{U_{\alpha \beta}}
 \isom 
(\mcP_{G, \epsilon} |_{U_\alpha}, \nabla_{\mcP_G, \epsilon} |_{U_\alpha}+ \epsilon \cdot  b_\alpha) |_{U_{\alpha \beta}} \notag
\end{align}
over $\tau^v_{C_S, \alpha \beta}$ for $(\alpha, \beta) \in I_2$.
The curve $C_{S}^v$ together with the resulting flat $G$-torsor, which we denote by $(\mcP_G^v, \nabla_{\mcP_G^v})$,
 specifies the deformation of $(C_S, \mcP_G, \nabla_{\mcP_G})$ 
 classified by $t_{C_S, \mcP_G, \nabla_{\mcP_G}}^{-1}(v)$.
 That is to say, 
 the assignment $v \mapsto (C_S^v, \mcP_G^v, \nabla_{\mcP^v_G})$
 gives the inverse of the bijection $t_{C_S, \mcP_G, \nabla_{\mcP_G}}$.

\vspace{3mm}

\subsection{} \label{Ws1599}

In this subsection, we discuss deformations  of indigenous bundles on a {\it Riemann surface}.
If,  in Definition \ref{D01} (i),  ``$C_S$"   is replaced by a Riemann surface and the words ``$B$-torsor" and ``connection" are understood in the analytic sense, then one obtains the notion of an  indigenous bundle on a Riemann surface.
Notice that the complex analytic stack $\mfS_{g, \mbC}^\text{an}$ associated with $\mfS_{g, \mbC}$  may  be thought of as the moduli stack  classifying connected compact Riemann surfaces of genus $g$ equipped with an  indigenous bundle.

Let $C_\mbC$ be a proper curve over $\mbC$ of genus $g$ and $\mcP^\circledast := (\mcP_B, \nabla_{\mcP_G})$ an indigenous bundle on $C_\mbC$.
 Then, $\mcP^\circledast$ gives rise to an indigenous bundle $\mcP^{\circledast,  \mr{an}} := (\mcP_B^\mr{an}, \nabla_{\mcP_G^\mr{an}})$ 
  on the Riemann surface  $C_\mbC^\mr{an}$  via the GAGA principle.
Just as in the case of the algebraic setting discussed before,
a $\mbC$-connection $\nabla^{\mr{ad}}_{\mcP_G^\mr{an}}$  
 on the (holomorphic) adjoint vector bundle $\mr{ad}(\mcP^\mr{an}_G)$ of the $G$-torsor $\mcP^\mr{an}_G$ are defined.
Also, 
we obtain a $\mbC$-linear morphism 
$\widetilde{\nabla}^\mr{ad}_{\mcP_B^\mr{an}} : \widetilde{\mcT}_{\mcP_B^\mr{an}} \migi \Omega_{C_\mbC}\otimes \mr{ad}(\mcP_G^\mr{an})$ and an isomorphism
of complexes 
$\eta^{\bullet \mr{an}} : \mcK^\bullet [\widetilde{\nabla}^\mr{ad}_{\mcP_B^\mr{an}}] \isom \mcK^\bullet [\nabla_{\mcP_G^\mr{an}}^\mr{ad}]$.

Next, denote by
 \begin{align}
 T_{C_\mbC^\mr{an}, \mcP^{\circledast, \mr{an}}} \notag
 \end{align}
the  deformation space of $(C_\mbC^\mr{an}, \mcP_B^\mr{an}, \nabla_{\mcP_G^\mr{an}})$ (as a data consisting of a Riemann surface $C_\mbC^\mr{an}$, a $B$-torsor over it, and a $\mbC$-connection on the $G$-torsor associated to this $B$-torsor) over $\mr{pt}_\epsilon := \mr{Spec}(\mbC [\epsilon]/(\epsilon^2))^\mr{an}$.
Then, 
it may be identified with the tangent space of $\mfS_{g, \mbC}^\mr{an}$ at the point classified by $(C_\mbC^\mr{an}, \mcP^{\circledast, \mr{an}})$. 
Also,  there exists a canonical isomorphism
\begin{align} \label{Eghkyi}
t_{C_\mbC^\mr{an}, \mcP^{\circledast, \mr{an}}}: T_{C_\mbC^\mr{an}, \mcP^{\circledast, \mr{an}}} \isom \mbH^1 (\mcK^\bullet [\widetilde{\nabla}^\mr{ad}_{\mcP_B^\mr{an}}]),
\end{align}
which makes the following diagram commute:
\begin{align} \label{WW1002}
\begin{CD}
T_{C_\mbC, \mcP^\circledast} @>t_{C_\mbC, \mcP^\circledast}>> \mbH^1 (\mcK^\bullet [\widetilde{\nabla}_{\mcP_B}^\mr{ad}]) @> \mbH^1 (\eta^{\bullet}) > \sim >  \mbH^1 (\mcK^\bullet [\nabla^\mr{ad}_{\mcP_G}])
\\
@VVV @VVV @VVV
\\
T_{C_\mbC^\mr{an}, \mcP^{\circledast, \mr{an}}} @> \sim > t_{C_\mbC^\mr{an}, \mcP^{\circledast, \mr{an}}} > \mbH^1 (\mcK^\bullet [\widetilde{\nabla}^\mr{ad}_{\mcP_B^\mr{an}}]) @> \sim > \mbH^1 (\eta^{\bullet, \mr{an}})  >  \mbH^1 (\mcK^\bullet [\nabla^\mr{ad}_{\mcP_G^\mr{an}}]), 
\end{CD} 
\end{align}
where
$\mbH^1 (\eta^{\bullet})$ denotes the isomorphism induced from $\eta^{\bullet}$ 
and  the vertical arrows are obtained naturally via the GAGA principle.
In fact, the explicit construction of this isomorphism can be given in the same way as discussed in the previous subsection. 
Since the rightmost  vertical arrow is an isomorphism (cf. ~\cite{De}, II, Th\'{e}or\`{e}me 6.13),
the remaining vertical arrows are in fact isomorphisms.

\vspace{10mm}
\section{Dormant Indigenous Bundles} \label{s16}\vspace{3mm}

In this section, we recall the definition of a {\it dormant} indigenous bundle
and discuss various moduli functors 
concerning dormant indigenous bundles.

\vspace{3mm}
\subsection{}

  Let  $S$, $G$, and $B$ 
  be as in \S\,\ref{s11}.
 Suppose further that $R = K$ for a field of characteristic $p>2$.
First, we recall the definition of $p$-curvature map.
Let 
$\pi : \mcP \migi C_S$ be a $G$-torsor over $C_S$ and
 $\nabla_\mcP : \mcT_{C_S/S} \migi \widetilde{\mcT}_{\mcP/S}$ 
an  $S$-connection  on $\mcP$.
If $\partial$ is a derivation corresponding to a local section $\partial$ of $\mcT_{C_S/S}$ 
(resp., $\widetilde{\mcT}_{\mcP/S}:= (\pi_*(\mcT_{\mcP/S}))^G)$,
then we shall denote by $\partial^{[p]}$ the $p$-th iterate of $\partial$, which specifies  a derivation corresponding to a local section of $\mcT_{C_S/S}$ (resp., $\widetilde{\mcT}_{\mcP/S}$).
Since the equality $\alpha_\mcP(\partial^{[p]})= (\alpha_\mcP(\partial))^{[p]}$ holds  for any local section $\partial$ of $\mcT_{C_S/S}$,
the image of the $p$-linear map $\mcT_{C_S/S} \migi \widetilde{\mcT}_{\mcP/S}$ defined by assigning $\partial \mapsto \nabla_\mcP (\partial^{[p]})-(\nabla_\mcP(\partial))^{[p]}$ is contained  in $\mr{ad}(\mcP)$ ($=\mr{Ker}(\alpha_\mcP))$.
The morphism of sheaves
\begin{align} 
  \psi_{\mcP, \nabla_\mcP} : \mcT_{C_S/S}^{\otimes p} \migi \mr{ad}(\mcP)  \notag
  \end{align}
determined by assigning
\begin{align} 
  \partial^{\otimes p} \mapsto  \nabla_\mcP (\partial^{[p]})-(\nabla_\mcP(\partial))^{[p]}. \notag
  \end{align}
is verified to be $\mcO_{C_S}$-linear.
We shall refer to this  morphism $\psi_{\mcP, \nabla_\mcP}$ as the {\bf $p$-curvature map} of $(\mcP,\nabla_\mcP)$.

\vspace{3mm}
\bde \label{d01}
 \begin{itemize}
 \item[(i)]
 We shall say that an indigenous bundle $\mcP^\ind := (\mcP_B,\nabla_{\mcP_G})$ on $C_S/S$
  is {\bf dormant} if the $p$-curvature map $\psi_{\mcP_G, \nabla_{\mcP_G}}$ of $(\mcP_G,\nabla_{\mcP_G})$ vanishes identically on $C_S$.
 \item[(ii)]
 Let $T$ be a $K$-scheme. A {\bf dormant curve}  over $T$ of genus $g$ is a pair 
 \begin{align}
 X_{/T}^{^\text{Zzz...}} := (X/T, \mcP^\ind) \notag
 \end{align}
   consisting of a proper  curve $X$ over $T$  of genus $g$ and a dormant indigenous bundle $\mcP^\ind$ on $X/T$.
 
 \item[(iii)]
 Let $T$ be a $K$-scheme,  and let $X_{/T}^{^\text{Zzz...}} := (X/T, \mcP^\ind := (\pi_\mcP : \mcP_B \migi X, \nabla_{\mcP_\mcG}))$ and  $Y_{/T}^{^\text{Zzz...}} := (Y/T, \mcQ^\ind :=(\pi_\mcQ : \mcQ_B\migi Y, \nabla_{\mcQ_\mcG}))$  be dormant curves  over $T$ of genus $g$.
 An {\bf isomorphism of dormant curves} from $X_{/T}^{^\text{Zzz...}} $ to $Y_{/T}^{^\text{Zzz...}}$ is a pair $(h, h_B)$ consisting of an isomorphism $h : X \isom Y$ of $T$-schemes and an isomorphism $h_B : \mcP_B \isom \mcQ_B$ that
 makes the square diagram 
  \begin{align} 
  \begin{CD}
 \mcP_B @> h_B > \sim > \mcQ_B
 \\
 @V\pi_\mcP VV @VV \pi_\mcQ V
 \\
 X @> \sim > h > Y
 \end{CD} \notag \end{align}
commute and 
  is compatible with both  the 
 $B$-actions and 
 the  $S$-connections in the evident sense.  
 \end{itemize}
 \ede

Next, we shall introduce the notion of   ordinariness for  dormant curves.
If $\mcP^\ind := (\mcP_B, \nabla_{\mcP_G})$ is an indigenous bundle on $C_S/S$,
i.e.,  ${C_S}^{^\text{Zzz...}}_{/S} :=(C_S/S, \mcP^\ind)$ form a dormant curve,  then
the natural morphism $\mr{Ker}(\nabla_{\mcP_G}^\mr{ad})[0] \migi \mcK^\bullet [\nabla_{\mcP_G}^\mr{ad}]$ determines, via the functor $\mbR^1 f_{S*} (-)$, 
 a morphism
 \begin{align} 
   \gamma^\natural_{\mcP^\ind} :\mbR^1f_{S*}(\mr{Ker}(\nabla_{\mcP_G}^\mr{ad})) \migi  \mbR^1f_{S*}(\mcK^\bullet [\nabla_{\mcP_G}^\mr{ad}])   \notag
   \end{align}
 of $\mcO_S$-modules.
By composing it and $\gamma_{\mcP^\ind}^\flat$ (cf. (\ref{WW1003})),
  we obtain a morphism
\begin{align} 
 \gamma_{\mcP^\ind}^\heartsuit  : \mbR^1f_{S*}(\mr{Ker}(\nabla_\mr{ad})) \migi \mbR^1f_{S*}(\mcT_{C_S/S})   \notag
 \end{align}
of $\mcO_S$-modules.
Notice  that $\gamma^\heartsuit_{\mcP^\ind}$ coincides with  the morphism obtained by applying the functor $\mbR^1f_{S*}(-)$ to the natural composite 
\begin{align} 
\mr{Ker}(\nabla_{\mcP_G}^\mr{ad}) \migiincl \mr{ad}(\mcP_G) \migisurj \mr{ad}(\mcP_G)/\mr{ad}(\mcP_G)^1 \xrightarrow{\overline{\nabla}_{\mcP_G}^\flat} \mcT_{C_S/S}. \notag
\end{align}

\vspace{3mm}
\bde  \label{d02}
 We shall say that  ${C_S}^{^\text{Zzz...}}_{/S} =(C_S/S, \mcP^\ind)$ is {\bf ordinary} if $\gamma_{\mcP^\ind}^\heartsuit$ is an isomorphism.
   \ede


\vspace{3mm}
\subsection{} \label{s18}

Denote by
\begin{align}  \label{e46}
\Dp_{g, K}  \  \ \ \left(\text{resp.}, {^\circledcirc \Dp_{g, K}}\right) 
\end{align}
the stack classifying dormant curves (resp., ordinary dormant curves) over $K$ of genus $g$. 
We obtain   a natural sequence of stacks
\begin{align}
 {^\circledcirc \Dp_{g, K}} \migi \Dp_{g, K} \migi   \mfS_{g, K}. \notag
 \end{align}
We here recall the following  result from the $p$-adic Teichm\"{u}ller theory studied by  S. Mochizuki.

\vspace{3mm}
\bpr \label{p05}
\begin{itemize}
\item[(i)]
  The stack  $\Dp_{g,K}$ may be represented by a nonempty, geometrically connected, and smooth Deligne-Mumford stack over $K$ of dimension $3g-3$, and forms  a  closed substack of $\mfS_{g,K}$.
 Moreover,  the projection $\Dp_{g,K} \migi \mfM_{g, K}$  is  finite and faithfully flat.
\item[(ii)]
 The stack  ${^\circledcirc \Dp_{g,K}}$ may be repressnted by  a dense  open  substack of $\Dp_{g,K}$ 
  and coincides with the \'{e}tale locus of
  $\Dp_{g,K}$  over $\mfM_{g,K}$. 
  In particular,  ${^\circledcirc \Dp_{g,K}}$ is a  nonempty, geometrically connected, and smooth  Deligne-Mumford stack over $K$ of dimension $3g-3$.
 \end{itemize}
  \epr
\begin{proof}
See
  ~\cite{Mzk2}, Chap.\,II, \S\,2.3, Lemma 2.7 and 
 Theorem 2.8 (and its proof). 
\end{proof}
\vspace{3mm}

Finally,  by means of cohomological expressions, we describe the differential of  the closed immersion   ${\Dp_{g, K}} \migiincl   \mfS_{g, K}$.
\vspace{3mm}
\bpr \label{p07} 
Let $\mcP^\ind$ be an indigenous bundle on $C_S/S$ and let $v : S \migi \mfS_{g, K}$ be (as in \S\,\ref{s15}) the $S$-rational point of $\mfS_{g, K}$ classifying $(C_S/S, \mcP^\ind)$.
Suppose further that 
$v$ factors through the closed immersion  ${\Dp_{g, K}} \migi   \mfS_{g, K}$, i.e., that $\mcP^\ind$ is dormant.
Denote by $\breve{v} : S \migi {\Dp_{g, K}}$ the resulting $S$-rational point of ${\Dp_{g, K}}$.
Then,  there exists a canonical isomorphism
\begin{align} 
      \breve{t}_{C_S, \mcP^\ind} :   \breve{v}^*(\mcT_{\Dp_{g,K}/K} ) \isom \mbR^1f_{S*}(\mr{Ker}(\nabla_{\mcP_G}^\mr{ad}))    \notag \end{align}
which makes the  following square diagram  commute:
\begin{align} 
\begin{CD}
 \breve{v}^*(\mcT_{\Dp_{g,K}/K} )  @>>>  v^*(\mcT_{\mfS_{g,K}/K})
\\
@V \breve{t}_{C_S, \mcP^\ind} V \wr V @V \wr V  V
\\
\mbR^1f_{S*}(\mr{Ker}(\nabla_{\mcP_G}^\mr{ad})) @> \gamma_{\mcP^\ind}^\natural >> \mbR^1f_{S*}(\mcK^\bullet [\nabla_{\mcP_G}^\mr{ad} ]), 
\end{CD} \notag \end{align}
where the upper horizontal arrow denotes the $\mcO_S$-linear  morphism obtained by differentiating 
the closed immersion   $\Dp_{g,K} \migiincl \mfS_{g, K}$ and the right-hand vertical arrow 
is the middle vertical arrow in (\ref{WW1004}).
  \epr
\begin{proof}
See the proof of ~\cite{Mzk2}, Chap.\,II, \S\,2.3, Theorem 2.8.
\end{proof}

\vspace{10mm}
\section{Comparison of Symplectic Structures} \label{s20}\vspace{3mm}

In this section, we construct a natural  symplectic structure (cf. (\ref{e5f5})) on the moduli stack $\mfS_{g, R}$ appearing in  the statement of  Theorem A (= Theorem \ref{t05}).

\vspace{3mm}
\subsection{} \label{s21}

Let  $R$ be as in \S\,\ref{s305},  and  $G$,  $B$  as in \S\,\ref{s11}.
Also, let $\overline{v} : S \migi \mfM_{g,R}$  be an object of $\mfS \mfc \mfh_{/\mfM_{g,R}}$ and 
 $\mcP^\ind :=(\mcP_B, \nabla_{\mcP_G})$ an indigenous bundle on $C_S/ S$.
Denote by 
 $v : S \migi \mfS_{g, R}$ 
the $S$-rational point of $\mfS_{g, R}$ classifying the pair $(C_S, \mcP^\ind)$.
Recall that the Killing form  on $\mfg$ ($= \mfs \mfl_2$)  is a {\it nondegenerate} symmetric bilinear map  $\kappa : \mfg \times \mfg \migi k$ defined by $\kappa (a, b) = \frac{1}{4} \cdot \mr{tr}(\mr{ad}(a)\cdot \mr{ad}(b))$ ($= \mr{tr}(ab)$) for any $a, b \in \mfg$, where we recall the assumption that $2$ is invertible in $R$.
The adjoint bundle $\mr{ad}(\mcP_G)$ has   an $\mcO_S$-bilinear  map
\begin{align} \label{e52}
 \kappa_{\mcP^\ind} : \mr{ad}(\mcP_G)\otimes \mr{ad}(\mcP_G) \migi \mcO_{C_S},
 \end{align}
obtained by  the  change of structure group via $\kappa$; it
 induces an isomorphism
\begin{align} 
  \kappa_{\mcP^\ind}^\vartriangleright : \mr{ad}(\mcP_G) \isom \mr{ad}(\mcP_G)^\vee. \notag
  \end{align}
Let us write 
  $ \nabla^{\mr{ad}, \otimes 2}_{\mcP_G}$ for the $S$-connection on the tensor product $\mr{ad}(\mcP_G)\otimes  \mr{ad}(\mcP_G)$ induced  naturally by $\nabla^\mr{ad}_{\mcP_G}$.
The  morphism $\kappa_{\mcP^\ind}$ is compatible with the respective $S$-connections  $ \nabla^{\mr{ad}, \otimes 2}_{\mcP_G}$ and $d$.
By composing  the morphism $\kappa_{\mcP^\ind}$ and the cup product in the de Rham cohomology, we obtain  a skew-symmetric $\mcO_S$-bilinear map on $\mbR^1f_{S*}(\mcK^\bullet [{\nabla_{\mcP_G}^\mr{ad}}])$:
\begin{align} \label{e54}
\mbR^1f_{S*}(\mcK^\bullet [{\nabla_{\mcP_G}^\mr{ad}}]) \otimes  
\mbR^1f_{S*}(\mcK^\bullet [{\nabla_{\mcP_G}^\mr{ad}}])
& \migi  \ \mbR^2f_{S*} (\mcK^\bullet [\nabla^{\mr{ad}, \otimes 2}_{\mcP_G}) \\
 & \migi  \ \mbR^2f_{S*}(\mcK^\bullet [d]) \notag  \\ 
 & \isom  \ \mbR^1f_{S*}(\Omega_{C_S/S}) \notag \\
 & \hspace{-0.5mm} \xrightarrow{\int_{C_S}}  \ \mcO_S \notag \end{align}
(cf.  ~\cite{ILL}, Corollary 5.6,  for the definition of  the third arrow).
The $\mcO_S$-bilinear maps (\ref{e54})
of the case where $(S, \mcP^\circledast)  = (\mfS_{g, R}, \mcP^\circledast_{g, R})$, where $\mcP^\circledast_{g, R}$ denotes the tautological indigenous bundle on the curve $C_{\mfS_{g, R}}$ over $\mfS_{g, R}$,  determines
  a $2$-form 
\begin{align} \label{e5f5}
 \omega^{\mr{PGL}}_{{g,R}}   
\end{align}
on 
$\mfS_{g,R}$ under the identification $\mcT_{\mfS_{g, R}} \isom \mbR^1 f_{\mfS_{g,R}*}(\mcK^\bullet [\nabla^\mr{ad}_{\mcP^\circledast_{g, R}}])$ resulting from 
 Proposition \ref{p02}.

\vspace{3mm}
\subsection{} \label{sy21}

In this subsection, we shall prove that $\omega^{\mr{PGL}}_{{g,R}}$ forms a symplectic structure on $\mfS_{g,R}$ (cf. Proposition \ref{p10}).
To this end, 
we first  prepare for describing Lemma \ref{pff27} below.

Let us consider the case where $R = \mbC$.
Let $C^\mr{an}$ be a connected compact Riemann surface of genus $g$,   $\Sigma$ its underlying   orientable closed  surface,   and $\pi := \pi_1(\Sigma, z)$ the fundamental group of $\Sigma$ with respect to a fixed base-point $z \in \Sigma$.
The space $\mr{Hom}(\pi, G)$ of representations $\pi \migi G$ has a canonical $G$-action  obtained by composing representations with inner automorphisms of $G$; we shall denote by $\mr{Hom}(\pi, G)/G$  the orbit space.
The subspace
\begin{align}
  \mfR  \ \left(\subseteq   \mr{Hom}(\pi, G)/G \right)\notag
\end{align}
 consisting of all irreducible representations forms a  complex manifold of dimension $6g-6$.
The Riemann-Hilbert correspondence gives a canonical identification between $\mfR$ with the moduli space of irreducible flat $G$-torsors (in the analytic sense) over $C^\mr{an}$ (cf.  ~\cite{Bi3}, \S\,2).
 By taking the monodromy representations   of indigenous bundles,
we obtain (cf. ~\cite{Hu}, Theorem) a local biholomorphic map 
\begin{align} 
  \mu :  \mfS_{g, \mbC}^\mr{an} \migi \mfR. \notag
  \end{align}

Now, let us take  a representation $\pi \migi G$, or equivalently,
an irreducible  flat $G$-torsor $(\mcP_G^\mr{an}, \nabla_{\mcP_G^\mr{an}})$ over  $C^\mr{an}$, 
 classified by a point of $\mfR$.
Denote by 
\begin{align}
T_{\mcP_G^\mr{an}, \nabla_{\mcP_G^\mr{an}}}
\end{align}
the tangent space of $\mfR$ at this point, i.e., the   deformation space   of $(\mcP_G^\mr{an}, \nabla_{\mcP_G^\mr{an}})$ over $C^\mr{an}_\epsilon \ \left(:= C^\mr{an} \times \mr{pt}_\epsilon \right)$.
Just as in the algebraic case (i.e., $t_{\mcP_G, \nabla_{\mcP_G}}$) discussed in \S\,\ref{Ws35},
there exists a canonical isomorphism
\begin{align}
t_{\mcP^\mr{an}_G, \nabla_{\mcP^\mr{an}_G}} : T_{\mcP_G^\mr{an}, \nabla_{\mcP_G^\mr{an}}} \isom \mbH^1 (\mcK^\bullet [\nabla^\mr{ad}_{\mcP^\mr{an}_G}]).
\end{align}
Then, we have the following lemma.

\vspace{3mm}
\ble \label{pff27}
The following square diagram is commutative:
\begin{align} \label{Ww470}
\begin{CD}
T_{C^\mr{an}, \mcP^{\circledast, \mr{an}}}@> t_{C^\mr{an}, \mcP^{\circledast, \mr{an}}} > \sim >
\mbH^1 (\mcK^\bullet [\widetilde{\nabla}^\mr{ad}_{\mcP_B^\mr{an}}])
\\
@V d\mu_{C^\mr{an}, \mcP^{\circledast, \mr{an}}} V \wr V @V \wr V \mbH^1 (\eta^{\bullet, \mr{an}})V
\\
T_{\mcP^\mr{an}_G, \nabla_{\mcP^\mr{an}_G}} @> \sim > t_{\mcP^\mr{an}_G, \nabla_{\mcP^\mr{an}_G}} > \mbH^1 (\mcK^\bullet [\nabla^\mr{ad}_{\mcP_G^\mr{an}}])
\end{CD}
\end{align} 
(cf. (\ref{Eghkyi}) for the definition of $t_{C^\mr{an}, \mcP^{\circledast, \mr{an}}}$), where the left-hand vertical arrow  $d\mu_{C^\mr{an}, \mcP^{\circledast, \mr{an}}}$ denotes the isomorphism  arising from the local biholomorphic map $\mu$.
\ele
\begin{proof}
Let us take an element $v$  of $\mbH^1 (\mcK^\bullet [\widetilde{\nabla}^\mr{ad}_{\mcP_B^\mr{an}}])$, which may be represented by  a $1$-cocycle $(\{ a_{\alpha, \beta} \}_{\alpha, \beta}, \{ b_\alpha \}_\alpha)$ 
 of the total complex $\mr{Tot}^\bullet (\check{C}^\bullet (\mfU, \mcK^\bullet [\widetilde{\nabla}^\mr{ad}_{\mcP_B^\mr{an}}]))$ (as in (\ref{W105})) with respect to a suitable analytic open covering $\mfU := \{ U_\alpha \}_{\alpha \in I}$ of $C^\mr{an}$.
Denote by 
$(C^{\mr{an}, v}, \mcP_{B}^{\mr{an}, v}, \nabla_{\mcP_{G}^{\mr{an}, v}})$ the deformation  
classified by the point of $T_{C^\mr{an}, \mcP^{\circledast, \mr{an}}}$ corresponding to $v$ via $t_{C^\mr{an}, \mcP^{\circledast, \mr{an}}}$.
In particular,  the underlying Riemann surface $C^{\mr{an}, v}$ of $(\mcP_{G}^{\mr{an}, v}, \nabla_{\mcP_{G}^{\mr{an}, v}})$, where $\mcP_{G}^{\mr{an}, v} := \mcP_{B}^{\mr{an}, v} \times^B G$, 
 corresponds to
the element of $H^1  (C^\mr{an}, \mcT_{C_\mbC^\mr{an}})$ represented by the $1$-cocycle $\{\alpha_{\mcP_G^\mr{an}}(a_{\alpha\beta}) \}_{\alpha, \beta}$.
We shall denote by $(\breve{\mcP}_{G}^{\mr{an}, v}, \breve{\nabla}_{\mcP_{G}^{\mr{an}, v}})$ the deformation of $(\mcP_G^\mr{an}, \nabla_{\mcP_G^\mr{an}})$ classified the image of
$v$ via $d\mu_{C_\mbC^\mr{an}, \mcP^{\circledast, \mr{an}}}$.
Then, after possibly replacing $\mfU$ with its refinement,  the element of 
$\mbH^1 (\mcK^\bullet [\nabla^\mr{ad}_{\mcP_G^\mr{an}}])$ corresponding to this deformation via $t_{\mcP^\mr{an}_G, \nabla_{\mcP^\mr{an}_G}}$  may be represented by a $1$-cocycle  $\breve{v} := (\{ \breve{a}_{\alpha \beta} \}_{\alpha, \beta}, \{ \breve{b}_\alpha \}_\alpha)$  of 
$\mr{Tot}^\bullet (\check{C}^\bullet (\mfU, \mcK^\bullet [\nabla^\mr{ad}_{\mcP_G^\mr{an}}]))$.
By the definition of $d\mu_{C_\mbC^\mr{an}, \mcP^{\circledast, \mr{an}}}$,  
   the pair  $(\breve{\mcP}_{G}^{\mr{an}, v}, \breve{\nabla}_{\mcP_{G}^{\mr{an}, v}})$
may be  characterized as a unique deformation of $(\mcP_G^\mr{an}, \nabla_{\mcP_G^\mr{an}})$ over $C^\mr{an}_\epsilon$ 
having  the same monodromy as $(\mcP_{G}^{\mr{an}, v}, \nabla_{\mcP_{G}^{\mr{an}, v}})$.
Hence, one verifies from the construction of the Riemann-Hibert correspondenc (cf, ~\cite{Chen}, Proposition 5.2.1) that  the equality $a_{\alpha \beta} - \nabla_{\mcP_G^\mr{an}}(\alpha_{\mcP_G^\mr{an}}(a_{\alpha\beta})) = \breve{a}_{\alpha \beta}$ holds for any $(\alpha, \beta) \in I_2$. 
 This implies  the equality $\mbH^1 (\eta^{\bullet, \mr{an}})(v)  = \breve{v}$ by the definition of $\eta^{\bullet, \mr{an}}$, and hence,
 completes the proof of the lemma.
\end{proof}
\vspace{3mm}

By means of the above lemma, one may prove the following proposition.

\vspace{3mm}
\bpr \label{p10}
$\omega^{\mr{PGL}}_{{g,R}}$ is 
a symplectic structure on $\mfS_{g, R}$. 
 \epr
\begin{proof}
First, let us prove that $\omega^{\mr{PGL}}_{{g,R}}$ is  nondegenerate.
Let $S$ be a scheme over $R$ and $v : S \migi \mfS_{g, R}$ an $S$-rational point, which classifies  a pair $(C_S, \mcP^\circledast)$ as before.
The isomorphisms $\int_{C_S, \mr{ad}(\mcP_G)}$ and   $\int_{C_S, \Omega_{C_S/S}\otimes\mr{ad}(\mcP_G)}$ (cf. (\ref{e16})) induces, under the identification 
  $\mr{ad}(\mcP_\mcG) = \mr{ad}(\mcP_G)^\vee$  arising from  $\kappa_{\mcP^\circledast}$,      isomorphisms
\begin{align}
 \mbR^qf_{S*}(\mcK^p [{\nabla_{\mcP_G}^\mr{ad}}]) \isom  \mbR^{1-q}f_{S*}(\mcK^{1-p} [{\nabla_{\mcP_G}^\mr{ad}}])^\vee
\end{align}
for the  pairs $(p, q)$ with $0 \leq p, q \leq 1$.
 The collection of these isomorphisms determines
 an isomorphism of spectral sequences from
``${E}_1^{p,q} = \mbR^qf_{S*}(\mcK^p [{\nabla_{\mcP_G}^\mr{ad}}]) \Rightarrow \mbR^{p+q}f_{S*}(\mcK^\bullet [{\nabla_{\mcP_G}^\mr{ad}}])$''
 to 
``${E}_1^{p,q} = \mbR^{1-q}f_{S*}(\mcK^{1-p} [{\nabla_{\mcP_G}^\mr{ad}}])^\vee \Rightarrow \mbR^{2-p-q}f_{S*}(\mcK^\bullet [{\nabla_{\mcP_G}^\mr{ad}}])^\vee$''.
Hence, it induces  an isomorphism $\mbR^1f_{S*}(\mcK^\bullet [{\nabla_{\mcP_G}^\mr{ad}}]) \isom  \mbR^1f_{S*}(\mcK^\bullet [{\nabla_{\mcP_G}^\mr{ad}}])^\vee$.
One verifies immediately that
the morphism $\mbR^1f_{S*}(\mcK^\bullet [{\nabla_{\mcP_G}^\mr{ad}}])^{\otimes 2} \migi \mcO_S$ associated with this isomorphism coincides with  (\ref{e54}).
This implies that $ \omega^{\mr{PGL}}_{{g,R}}$ is nondegenerate.


Next, we consider the closedness of $\omega^{\mr{PGL}}_{{g,R}}$.
Since the closedness of a differential form is preserved under base-change via $\mr{Spec}(R) \migi \mr{Spec}(\mbZ[\frac{1}{2}])$, it suffices to verify the case where $R = \mbZ[\frac{1}{2}]$.
But, $\bigwedge^3 \Omega_{\mfS_{g, \mbZ[\frac{1}{2}]}/\mbZ[\frac{1}{2}]}$ is flat over $\mbZ[\frac{1}{2}]$, so the assertion of the case where $R = \mbZ [\frac{1}{2}]$ follows from that of the case where $R = \mbC$.
According to the discussion in  ~\cite{Go}, \S\,1.4,  \S\,1.7,  and Theorem (or ~\cite{Bi3}, \S\,3.2), $\mfR$ admits a  holomorphic symplectic structure 
\begin{align}
\omega_\mfR \in \Gamma (\mfR, {\bigwedge}^2 \Omega_{\mfR}).\notag
\end{align}
It follows from ~\cite{Bi3}, Theorem 3.2, that,  if  $\phi$ is a point of $\mfR$ and   $(\mcP_G^\mr{an}, \nabla_{\mcP_G^\mr{an}})$ denotes  the corresponding  flat $G$-torsor,  then 
the pairing
\begin{align}
\mbH^1 (\mcK^\bullet [\nabla^\mr{ad}_{\mcP_G^\mr{an}}]) \times
\mbH^1 (\mcK^\bullet [\nabla^\mr{ad}_{\mcP_G^\mr{an}}]) \migi \mbC \notag
\end{align}
 determined by  $\omega_\mfR$ under the isomorphism $t_{\mcP^\mr{an}_G, \nabla_{\mcP^\mr{an}_G}}$ may be obtained (up to a constant factor) in the same manner as (\ref{e54}).
In particular, the pull-back $\mu^*(\omega_\mfR)$ defines a symplectic structure on $\mfS_{g, \mbC}^\mr{an}$.
By Lemma \ref{pff27} together with the commutativity of the diagram (\ref{WW1002}), 
the $2$-form on $\mfS_{g, \mbC}$ corresponding to $\mu^*(\omega_\mfR)$
via the GAGA principle coincides with $\omega^\mr{PGL}_{g, \mbC}$.
In particular,  $\omega^\mr{PGL}_{g, \mbC}$ turns out to be closed.
This completes the proof of the proposition.
\end{proof}
\vspace{3mm}

\subsection{} \label{s22}

Let us keep the notation in \S\,\ref{s21}.
Suppose further that $\overline{v} : S \migi \mfM_{g,R}$ is \'{e}tale, hence $\mfS_S$ is a  smooth Deligne-Mumford stack over $R$ of relative dimension $6g-6$.
Consider the composite isomorphism
\begin{align}  \label{eE19}
\Omega_{S/R} \ \left(=\mcT_{S/R}^\vee \right)
& \isom \overline{v}^*(\mcT_{\mfM_{g,R}/R}^\vee) \\
& \isom \overline{v}^*(\mbR^1f_{g,R*}(\mcT_{C_{g,R}/\mfM_{g,R}})^\vee) \notag \\
&  \isom \overline{v}^*(f_{g,R*}(\Omega_{C_{g,R}/\mfM_{g,R}}^{\otimes 2})) \notag \\
& \isom f_{S*}(\Omega_{C_S/S}^{\otimes 2}), \notag
\end{align}
where
\begin{itemize}
\item[$\bullet$]
the first isomorphism follows from the \'{e}taleness of $\overline{v}$;
\item[$\bullet$]
the second and third isomorphisms
arises from 
the Kodaira-Spencer map of $C_{g, R}/\mfM_{g, R}$ (i.e., the rightmost  vertical arrow in (\ref{WW1004})) and 
$\int_{C_{g,R}, \Omega^{\otimes 2}_{C_{g,R}/\mfM_{g,R}}}$  (cf.  (\ref{e16})) respectively.

\end{itemize}
By this composite isomorphism,  the cotangent bundle $T^\vee_S$ ($= \mbA (\Omega_{S/R})$) may be thought of as the  {\it trivial} $\mbA(f_{S *}(\Omega_{C_S/S}^{\otimes 2}))$-torsor over $S$. 
Hence, for each global section $\sigma : S \migi \mfS_S$ of the  projection $\pi^\mfS_S : \mfS_S \migi S$,
there exists (cf. Proposition \ref{p01}) a unique isomorphism
\begin{align} \label{e58}
  \Psi_\sigma : T^\vee_S \isom \mfS_S 
  \end{align}
over $S$ which  is compatible with the respective $\mbA (f_{S*}(\Omega_{C_S/S}^{\otimes 2}))$-actions  and whose restriction $\Psi_\sigma |_{0_S}$ to the zero section $0_S : S \migi T^\vee_S$ coincides with $\sigma$.
In particular, $\Psi_\sigma$ induces an isomorphism of short exact sequences
\begin{align} 
 \begin{CD}
0 @>>> 0_S^*(\mcT_{T^\vee_S/S}) @>>> 0_S^*(\mcT_{T^\vee_S/R}) @>>> \mcT_{S/R} @>>> 0
\\
@. @V \wr VV @V\wr VV @V V  \mr{id} V @.
\\
0 @>>> \sigma^*(\mcT_{\mfS_S /S}) @>>> \sigma^*(\mcT_{\mfS_S/R}) @>>> \mcT_{S/R} @>>> 0.
\end{CD} \notag\end{align}

Next,  write $\overline{v}^\mfS :\mfS_S \migi \mfS_{g,R}$ for the base-change of $\overline{v}$ via the projection $\pi^\mfS_{g, R} : \mfS_{g,R} \migi \mfM_{g,R}$.
Since $\overline{v}^\mfS$ is \'{e}tale, which implies that   $\overline{v}^{\mfS*}(\Omega_{\mfS_{g,R}/R}) \cong \Omega_{\mfS_S/R}$,
 the pull-back
\begin{align} \label{e60}
 \omega_S^{\mr{PGL}} :=\overline{v}^{\mfS*}(\omega^{\mr{PGL}}_{{g,R}})
 \end{align}
 determines a symplectic structure on $\mfS_S$ (cf. Proposition \ref{p10}).

\vspace{3mm}
\subsection{} \label{s25}

Let $K$ be a field of characteristic $p >2$.
We shall write
\begin{align} \label{e100}
{^\circledcirc T^{\vee ^\mr{Zzz...}}_{g,K}} &  := T^\vee_{{^\circledcirc \Dp_{g,K}}} \ (= T^\vee_{\mfM_{g, K}} \times_{\mfM_{g, K}} {^\circledcirc \Dp_{g,K}}), \\
  {^\circledcirc \mfS^{^\mr{Zzz...}}_{g, K}} & := \mfS_{g,K} \times_{\mfM_{g,K}} {^\circledcirc \mfM^{^\mr{Zzz...}}_{g, K}}. \notag
\end{align}
The  stack  ${^\circledcirc \mfS^{^\mr{Zzz...}}_{g, K}}$ admits a section 
$ \sigma_{g,K} : {^\circledcirc \Dp_{g,K}} \migi {^\circledcirc \mfS^{^\mr{Zzz...}}_{g, K}}$ arising from the natural  immersion ${^\circledcirc \Dp_{g,K}} \migi \mfS_{g, K}$.
By Proposition \ref{p05},
one may apply the discussion in the previous subsection  to the case where 
the data ``$(R, S, \sigma : S \migi \mfS_S)$" is taken to be
$(K, {^\circledcirc \mfM^{^\mr{Zzz...}}_{g, K}}, \sigma_{g,K})$.
Thus,
we obtain an isomorphism
\begin{align} \label{e61}
 \Psi_{g,K} :   {^\circledcirc T^{\vee ^\mr{Zzz...}}_{g,K}} \isom  {^\circledcirc \mfS^{^\mr{Zzz...}}_{g, K}} 
 \end{align}
 over ${^\circledcirc \Dp_{g,K}}$ (cf. (\ref{e58})) and 
a symplectic structure
\begin{align} 
  \omega^{\mr{PGL}}_\circledcirc \notag
  \end{align}
 on ${^\circledcirc \mfS^{^\mr{Zzz...}}_{g, K}}$ (cf. (\ref{e60})).
On the other hand, 
 recall that ${^\circledcirc T^{\vee ^\mr{Zzz...}}_{g,K}}$ admits a canonical symplectic structure 
\begin{align} 
\omega^\text{can}_\circledcirc := \omega^\text{can}_{ {^\circledcirc \Dp_{g,K}}}\notag
\end{align}
 (cf. (\ref{e10})).
 The main assertion of the present paper, i.e.,  Theorem A, is  described as follows.
\vspace{3mm}
\bt \label{t05}
  The isomorphism $\Psi_{g,K}$ preserves  the  symplectic structure, i.e., the following equality holds: 
  \begin{align} 
  \Psi_{g, K}^*( \omega^{\mr{PGL}}_{\circledcirc} ) =  \omega^{\mr{can}}_{\circledcirc}.\notag
  \end{align}
\et

\vspace{3mm}
\begin{rema} \label{R071w}
 Let $X$ be a smooth Deligne-Mumford stack over $K$ equipped with a symplectic structure $\omega$.
 We shall say that a smooth substack $Y$ of $X$  is
 {\bf Lagrangian} (with respect to $\omega$) if $\omega |_Y =0$ and $\mr{dim}(Y) = \frac{1}{2} \cdot  \mr{dim}(X)$. 
 For example, 
 if $X = T_Y^\vee$, where $Y$ denotes a smooth  Deligne-Mumford stack over $K$ considered as a closed substack of $X$ by 
  the zero section $0_Y : Y \migi  T_Y^\vee \ \left(=X \right)$,
 then  $Y$ is Lagrangian with respect to the symplectic structure $\omega_Y^\mr{can}$.
  This example together with the  above theorem shows that  the substack ${^\circledcirc\mfM}_{g, K}^{^\mr{Zz...}}$ of $\mfS_{g, K}$ is Lagrangian with respect to the symplectic structure $\omega^\mr{PLG}_{g, K}$.
\end{rema}

\vspace{10mm}
\section{Proof of Theorem A} \label{s30}\vspace{3mm}

This section is devoted to  the proof of  Theorem \ref{t05}.
\vspace{3mm}
\subsection{} \label{s031}
To begin with,  we shall  discuss  the  translations of  the symplectic structures 
 $\omega^{\mr{PGL}}_{\circledcirc}$, $\omega^{\mr{can}}_{\circledcirc}$ with respect to  the affine structures on the respective underlying  spaces.
Let  $R$ be as in \S\,\ref{s305} and 
$\overline{v} : S \migi \mfM_{g,R}$
 an \'{e}tale relative scheme over $\mfM_{g, R}$.
In the following,  $(-)$ denotes either $T^\vee$ or $\mfS$.
We shall write
\begin{align}
(-)_S^{T^\vee} := (-)_S \times_S T^\vee_S,\notag
\end{align}
 which has
the projections 
\begin{align} 
\pi^{(-)}_{S, 1} : (-)_S^{T^\vee}  \migi (-)_S \ \  \text{and} \ \ 
\pi^{(-)}_{S, 2} : (-)_S^{T^\vee}  \migi T^\vee_S\notag
\end{align}
 onto the first  and second  factors respectively.
Define $\tau^{(-)}_{S}$ to be the composite 
\begin{align}
\tau_S^{(-)} : (-)_{S}^{T^\vee} 
 \migi 
  (-)_{S}^{T^\vee}  \times_{S} \mbA (\Omega_{S/R}) 
  \isom 
  (-)_{S}^{T^\vee}   \times_S \mbA (f_{S*} (\Omega_{C_S/S}^{\otimes 2}))
       \migi   (-)_{S}^{T^\vee}, \notag
\end{align}
where
\begin{itemize}
\item[$\bullet$]
the first arrow denotes  the product of the identity morphism of $(-)_S$ and   the section $T^\vee_S \migi T^\vee_S \times_S \mbA (\Omega_{S/R})$
 corresponding to  the Liouville form $\lambda_S \in \Gamma (T^\vee_S, \pi_S^{T^\vee *}(\Omega_{S/R}))$
(cf. \S\,\ref{s07}) on $T^\vee_S$;
\item[$\bullet$]
the second arrow  denotes  the product of the identity morphism
  of $(-)_{S}^{T^\vee}$ and the  isomorphism $\mbA (\Omega_{S/R}) \isom \mbA (f_{S*} (\Omega_{C_S/S}^{\otimes 2}))$ (cf. (\ref{eE19}));
\item[$\bullet$]
the third arrow
arises from the structure of $\mbA (f_{S *} (\Omega_{C_S/S}^{\otimes 2}))$-torsor (cf.  (\ref{eE19}) and Proposition \ref{p01}) on the first factors  in 
$(-)_S^{T^\vee} := (-)_S \times_S T^\vee_S$.
\end{itemize}
For each  $A  \in \Gamma (S, \Omega_{S/R})$ ($= \Gamma (C_S, \Omega_{C_S/S}^{\otimes 2})$),
we  denote by 
\begin{align} 
 \tau_{S, A}^{(-)} :  (-)_S \isom (-)_S\notag
 \end{align}
 the automorphism   of $(-)_S$
  determined, via its own  affine structure,  by the translation 
by $A$.
Also, denote by 
\begin{align} 
\sigma^{T^\vee}_{S, A} : S \migi T^\vee_S \notag
\end{align}
 the section of $\pi^{T^\vee}_S$ corresponding to $A$.
 In particular,     $0_S = \sigma_{S, 0}^{T^\vee}$.
By the various definitions involved, the equality
\begin{align} \label{eE68}
\tau^{(-)}_{S, A} = \pi_{S, 1}^{(-)} \circ \tau_S^{(-)} \circ (\mr{id}_{(-)_S} \times \sigma_{S, A}^{T^\vee})
\end{align}
 holds.
By replacing $S$ with 
  a covering space $S'$ of $\mfM_{g, \mbC}^\mr{an}$  (e.g., $\mfT^\Sigma$),
 we obtain, in the same manner as above,
 a complex analytic stack $(-)_{S'}^{T^\vee}$ and 
 various morphisms $\pi_{S', 1}^{(-)}$, $\pi_{S', 2}^{(-)}$, and $\tau_{S'}^{(-)}$.
 (They   will be used in the proof of Proposition  \ref{pf2f7}.)
 
 Then, the following Propositions \ref{pf2hf7} and \ref{pf2f7} hold.

\vspace{3mm}
\bpr \label{pf2hf7} 
 The following  equality  in $\Gamma (T^{\vee T^\vee}_{S}, \bigwedge^2 \Omega_{T^{\vee T^\vee}_{S}/R})$ holds:
\begin{align} \label{eE61}
 (\pi_{S, 1}^{T^\vee} \circ \tau_S^{T^\vee})^*(\omega_{S}^{\mr{can}})  = \pi_{S, 1}^{T^\vee *} (\omega_{S}^\mr{can}) + \pi^{T^\vee *}_{S, 2} (\omega^\mr{can}_S).
\end{align}
In particular, for each $A \in \Gamma (S, \Omega_{S/R})$, 
the following equality holds:
\begin{align} 
\tau^{T^\vee *}_{S, A} (\omega^\mr{can}_S) = \omega^\mr{can}_S + \pi_S^{T^\vee *}(dA).\notag
\end{align}
\epr
\begin{proof}
Let us  consider the former assertion.
Let $q_1, \cdots, q_{3g-3}$ be  local coordinates in $S$, and let  $p_1, \cdots, p_{3g-3}$  and $p'_1, \cdots p'_{3g-3}$ be  its  dual coordinates in the first and second  factors, respectively,  
of  the product $T^{\vee T^\vee}_S = T^\vee_S \times_S T^\vee_S$.
Then, locally on $T^{\vee T^\vee}_S$,  we have
\begin{align}
 (\pi_{S, 1}^{T^\vee} \circ \tau_S^{T^\vee})^*(\omega_{S}^{\mr{can}})  & =
 (\pi_{S, 1}^{T^\vee} \circ \tau_S^{T^\vee})^*(\sum_{i=1}^{3g-3} d p_i \wedge dq_i) \notag \\
 & =  \sum_{i=1}^{3g-3} d (p_i + p'_i) \wedge dq_i \notag   \\
&  = \left( \sum_{i=1}^{3g-3} d p_i \wedge dq_i \right)+ \left( \sum_{i=1}^{3g-3} d p'_i  \wedge dq_i \right) \notag \\
& =   \pi_{S, 1}^{T^\vee *} (\omega_{S}^\mr{can}) + \pi^{T^\vee *}_{S, 2} (\omega_S^\mr{can}). \notag
\end{align}
This completes the proof of the former assertion.

The latter assertion follows from the following sequence of equalities:
\begin{align} 
&  \ \ \ \ \,  \tau^{T^\vee *}_{S, A} (\omega^\mr{can}_S) \notag \\
& =  (\mr{id}_{T^\vee_S} \times \sigma_{S, A}^{T^\vee})^* ((\pi_{S, 1}^{T^\vee} \circ \tau_S^{T^\vee})^*(\omega^\mr{can}_S)) \notag  \\
& =  (\mr{id}_{T^\vee_S} \times \sigma_{S, A}^{T^\vee})^* (\pi_{S, 1}^{T^\vee *} (\omega_{S}^\mr{can}) +  \pi^{T^\vee *}_{S, 2} (\omega^\mr{can}_S)) \notag \\
& = (\pi_{S, 1}^{T^\vee}  \circ (\mr{id}_{T^\vee_S} \times \sigma_{S, A}^{T^\vee}))^* (\omega^\mr{can}_S) +  (\sigma_{S, A}^{T^\vee} \circ \pi^{T^\vee}_S)^* (d \lambda_S) \notag \\
& = \omega^\mr{can}_S + \pi^{T^\vee *}_S (d (\sigma^{T^\vee *}_{S, A}( \lambda_S))) \notag \\
& =\omega^\mr{can}_S + \pi_S^{T^\vee *}(dA), \notag 
\end{align}
where the first equality follows from (\ref{eE68}) and  the second equality follows from the former assertion, i.e., the equality  (\ref{eE61}).
\end{proof}
\vspace{3mm}

\bpr \label{pf2f7}
 The following equality  in $\Gamma (\mfS^{T^\vee}_S, \bigwedge^2 \Omega_{\mfS^{T^\vee}_S/R})$ holds:
 \begin{align} \label{E70002}
 (\pi_{S, 1}^{\mfS} \circ \tau_S^\mfS)^*(\omega_{S}^\mr{PGL})  = \pi^{\mfS*}_{S, 1} (\omega_{S}^\mr{PGL}) + \pi^{\mfS*}_{S, 2} (\omega^\mr{can}_S).
 \end{align}
 In particular, for each  $A \in \Gamma (S, \Omega_{S/R})$, 
 the following equality holds:
\begin{align} 
\tau_{S, A}^{\mfS*}(\omega^{\mr{PGL}}_{S}) = \omega^{\mr{PGL}}_{S} +\pi^{\mfS *}_S (dA).\notag
\end{align} 
\epr
\begin{proof}
We shall  prove the former assertion.
The equality (\ref{E70002}) may be obtained as the  pull-back, via the composite $S \stackrel{\overline{v}}{\migi} \mfM_{g, R} \migi \mfM_{g, \mbZ [\frac{1}{2}]}$,  of  the equality (\ref{E70002}) in the case  where  $(R, S)$ is taken to be $(\mbZ [\frac{1}{2}], \mfM_{g, \mbZ [\frac{1}{2}]})$.
Hence, 
the proof  is  reduced to 
the case where
$(R, S) = (\mbZ [\frac{1}{2}], \mfM_{g, \mbZ [\frac{1}{2}]})$.
Notice here that   the morphism 
\begin{align}
\Gamma (\mfS^{T^\vee}_{g, \mbZ [\frac{1}{2}]},  {\bigwedge}^2 \Omega_{\mfS^{T^\vee}_{g, \mbZ [\frac{1}{2}]}/\mbZ [\frac{1}{2}]}) \migi \Gamma (\mfS^{T^\vee}_{g,  \mbC},  {\bigwedge}^2 \Omega_{\mfS^{T^\vee}_{g, \mbC}/\mbC}),\notag
\end{align}
where $\mfS^{T^\vee}_{g, \mbZ [\frac{1}{2}]} := \mfS^{T^\vee}_{\mfM_{g, \mbZ [\frac{1}{2}]}}$, $\mfS^{T^\vee}_{g, \mbC} := \mfS^{T^\vee}_{\mfM_{g, \mbC}}$, arising from base-change via $\mr{Spec} (\mbC) \migi \mr{Spec} (\mbZ [\frac{1}{2}])$ is injective because $\bigwedge^2 \Omega_{\mfS^{T^\vee}_{g, \mbZ [\frac{1}{2}]}/\mbZ [\frac{1}{2}]}$ is flat over $\mbZ [\frac{1}{2}]$.
Thus,  it suffices
to prove the equality (\ref{E70002}) of  the case where 
$(R, S) = (\mbC, \mfM_{g, \mbC})$.
Moreover, by applying the analytification functor and replacing $\mfM_{g, \mbC}^\mr{an}$ by the Teichm\"{u}ller space $\mfT^\Sigma$, i.e., its universal covering,
one can reduce the problem to proving the equality 
 \begin{align} \label{Ef70002d}
 (\pi_{\mfT^\Sigma, 1}^{\mfS} \circ \tau_{\mfT^\Sigma}^{\mfS})^*(\omega_{\mfT^\Sigma}^{\mr{PGL}})  = \pi^{\mfS*}_{\mfT^\Sigma, 1} (\omega_{\mfT^\Sigma}^{\mr{PGL}}) + \pi^{\mfS*}_{\mfT^\Sigma, 2} (\omega^\mr{can}_{\mfT^\Sigma})
 \end{align}
 of holomorphic $2$-forms  on $\mfT^\Sigma$.

Since $\mfT^\Sigma$ is simply connected, 
there exists  a complex differential $1$-form $\delta$ on $\mfT^\Sigma$ satisfying the equality   $(1+ \sqrt{-1}) \cdot \omega_{WP} = - d \delta$,
  where $\omega_{WP}$ denotes the Weil-Petersson  symplectic form on $\mfT^\Sigma$.
Denote by $\sigma^F_{\mr{unif}}$ and $\sigma^B_\mr{unif}$ the 
sections 
$\mfT^\Sigma \migi \mfS_{\mfT^\Sigma}$
 determined by the Fuchsian and Bers uniformizations respectively   (cf. ~\cite{Ea}; ~\cite{AGBi}, \S\,3.1, (3.1)).
Because of the $\mbA (\Omega_{\mfT^\Sigma})$-torsor   structure on  $\mfS_{\mfT^\Sigma}$, it makes sense to speak of the sum  
\begin{align}
\sigma_{\mr{unif}, \delta}^B := \sigma^B_\mr{unif}  + \delta,\notag
\end{align}
  which specifies  a section $\mfT^\Sigma \migi \mfS_{\mfT^\Sigma}$.
The differences $\sigma^F_{\mr{unif}}  -\sigma^B_\mr{unif}$ and $\sigma^F_{\mr{unif}}  - \sigma_{\mr{unif}, \delta}^B$ may be thought of as elements of $\Gamma (\mfT^\Sigma, \Omega_{\mfT^\Sigma})$.
Then,
\begin{align}
d (\sigma^F_{\mr{unif}}  - \sigma_{\mr{unif}, \delta}^B) & = d (\sigma^F_{\mr{unif}}  -\sigma^B_\mr{unif}) - d \delta
  =  - \sqrt{-1} \cdot \omega_{WP}  + (1+ \sqrt{-1}) \cdot \omega_{WP} 
   =  \omega_{WP}, \notag
\end{align}
where  the second   equality follows from the definition of $\delta$ and ~\cite{McM}, Theorem 1.5.
Thus, it follows from ~\cite{Los}, Theorem 6.8,   that
the diffeomorphism $\Psi_{\sigma^B_{\mr{unif}, \delta}} : T^\vee_{\mfT^\Sigma} \isom \mfS_{\mfT^\Sigma}$ induced by $\sigma^B_{\mr{unif}, \delta}$ 
satisfies the equality $\Psi_{\sigma^B_{\mr{unif}, \delta}}^* (\omega^\mr{PGL}_{\mfT^\Sigma}) = \omega^\mr{can}_{\mfT^\Sigma}$.
Hence, the following sequence of equalities holds:
\begin{align} \label{eEE2}
 & \ \ \ \ (\Psi_{\sigma_{\mr{unif}, \delta}^B} \times \mr{id}_{T^\vee_{\mfT^\Sigma}})^*((\pi_{\mfT^\Sigma, 1}^{\mfS} \circ \tau^{\mfS}_{\mfT^\Sigma})^*(\omega_{\mfT^\Sigma}^\mr{PGL})) \\
&  =
( \pi_{\mfT^\Sigma, 1}^{\mfS} \circ \tau^{\mfS}_{\mfT^\Sigma} \circ (\Psi_{\sigma_{\mr{unif}, \delta}^B} \times \mr{id}_{T^\vee_{\mfT^\Sigma}}))^* (\omega_{\mfT^\Sigma}^\mr{PGL}) \notag \\
& = ( \pi_{\mfT^\Sigma, 1}^{\mfS} \circ (\Psi_{\sigma_{\mr{unif}, \delta}^B} \times \mr{id}_{T^\vee_{\mfT^\Sigma}}) \circ \tau^{T^\vee}_{\mfT^\Sigma})^* (\omega_{\mfT^\Sigma}^\mr{PGL}) \notag  \\
& = (\Psi_{\sigma_{\mr{unif}, \delta}^B} \circ \pi_{\mfT^\Sigma, 1}^{\mfS,  \mr{an}} \circ \tau^{T^\vee}_{\mfT^\Sigma})^* (\omega_{\mfT^\Sigma}^\mr{PGL}) \notag \\
& = ( \pi_{\mfT^\Sigma, 1}^{\mfS} \circ \tau^{T^\vee}_{\mfT^\Sigma})^*(\omega_{\mfT^\Sigma}^\mr{can}). \notag
\end{align}
On the other hand, we have the following equalities: 
\begin{align} \label{eEE1}
& \ \ \ \ (\Psi_{\sigma_{\mr{unif}, \delta}^B} \times \mr{id}_{T^\vee_{\mfT^\Sigma}})^* (\pi^{\mfS*}_{\mfT^\Sigma, 1} (\omega_{\mfT^\Sigma}^\mr{PGL}) + \pi^{\mfS*}_{\mfT^\Sigma, 2} (\omega^\mr{can}_{\mfT^\Sigma})) \\
& = (\pi^{\mfS*}_{\mfT^\Sigma, 1} \circ (\Psi_{\sigma_{\mr{unif}, \delta}^B} \times \mr{id}_{T^\vee_{\mfT^\Sigma}}))^*(\omega_{\mfT^\Sigma}^\mr{PGL}) + (\pi^{\mfS*}_{\mfT^\Sigma, 2} \circ (\Psi_{\sigma_{\mr{unif}, \delta}^B} \times \mr{id}_{T^\vee_{\mfT^\Sigma}}))^*(\omega^\mr{can}_{\mfT^\Sigma}) \notag  \\
& = (\Psi_{\sigma_{\mr{unif}, \delta}^B} \circ \pi^{T^\vee}_{\mfT^\Sigma, 1})^* (\omega_{\mfT^\Sigma}^\mr{PGL}) + \pi^{T^\vee*}_{\mfT^\Sigma, 2} (\omega^\mr{can}_{\mfT^\Sigma}) \notag \\
& = \pi^{T^\vee*}_{\mfT^\Sigma, 1} (\omega_{\mfT^\Sigma}^\mr{can}) + \pi^{T^\vee*}_{\mfT^\Sigma, 2} (\omega^\mr{can}_{\mfT^\Sigma}). \notag
\end{align}
By (\ref{eEE2}) and (\ref{eEE1}), the equality (\ref{Ef70002d}) holds if and only if the equality
\begin{align} \label{eEE3}
( \pi_{\mfT^\Sigma, 1}^{T^\vee} \circ \tau^{T^\vee}_{\mfT^\Sigma})^*(\omega_{\mfT^\Sigma}^\mr{can}) = \pi^{T^\vee*}_{\mfT^\Sigma, 1} (\omega_{\mfT^\Sigma}^\mr{can}) + \pi^{T^\vee*}_{\mfT^\Sigma, 2} (\omega^\mr{can}_{\mfT^\Sigma})
\end{align}
holds.
But, the equality (\ref{eEE3}) follows from an argument similar to the argument in the proof of the former assertion of Proposition \ref{pf2hf7}.
This completes the proof of the former assertion.

The latter assertion follows from the former assertion and  the  argument in the proof of the latter assertion of  Proposition \ref{pf2hf7}  where ``$T^\vee$" is replaced by ``$\mfS$". 
This completes the proof of the proposition.
\end{proof}

\vspace{3mm}
\subsection{} \label{s31}
Now, let us consider the case where  $R =K$ for a field $K$ of characteristic $p >2$.
Suppose further that  $\overline{v} : S \migi \mfM_{g, K}$ factors through the projection ${^\circledcirc\Dp_{g, K}} \migi \mfM_{g, K}$ and that the resulting morphism $\breve{v} : S \migi {^\circledcirc\Dp_{g, K}}$ is \'{e}tale.
Denote by
\begin{align} 
   \sigma_S : S \migi \mfS_S \ \ \left(\text{resp.,} \  \Psi_S :  T^\vee_S \migi  \mfS_S \right)  \notag\end{align}
the restriction of $\sigma_{g, K}$ (resp., $\Psi_{g, K}$) to $S$.
Since the natural map 
\begin{align} 
\Gamma ({^\circledcirc T_{g, K}^{\vee ^\mr{Zzz...}}},{\bigwedge}^2\Omega_{{^\circledcirc T_{g, K}^{\vee ^\mr{Zzz...}}}/K}) \migi \Gamma (T^\vee_S, {\bigwedge}^2\Omega_{T^\vee_S/K})\notag
\end{align}
induced from $\breve{v}$ is injective, the proof of Theorem \ref{t05} is reduced to proving the equality
\begin{align}
\Psi^*_S(\omega^\mr{PGL}_S) = \omega^\text{can}_S.\notag
\end{align}
Moreover, for the same reason,  we are always free to replace $S$ by any \'{e}tale covering of $S$. 

\vspace{3mm}
\subsection{} \label{s32}
It follows from  Propositions \ref{pf2hf7} and \ref{pf2f7}  that 
 for each  $A \in \Gamma (S, \Omega_{S/K})$,  we have   the following sequence of  equalities:
\begin{align} \label{E70007}
&  \ \ \ \  \sigma_S^*(\omega^\mr{PGL}_S)  - 0_S^*(\omega^{\mr{can}}_{S}) \\
 &  = \sigma^*_S (\tau_{S, A}^{\mfS *} (\omega_S^\mr{PGL}) - dA) - 0^*_S (\tau^{T^\vee}_{S, A} (\omega_S^\mr{can}) -dA) \notag  \\
 & = (0_S^* (\Psi_S^*(\tau_{S, A}^{\mfS *}(\omega^{\mr{PGL}}_{S}))) -dA) -
 (0^*_S (\tau^{T^\vee}_{S, A} (\omega_S^\mr{can})) -dA)
  \notag \\
  & = 0_S^* (\tau_{S, A}^{T^\vee*}(\Psi_S^*(\omega^{\mr{PGL}}_{S}))  -\tau^{T^\vee*}_{S, A}(\omega^{\mr{can}}_{S})) \notag \\
 & = \sigma_{S, A}^{T^\vee *} (\Psi_S^*(\omega^{\mr{PGL}}_{S}) - \omega^{\mr{can}}_{S}). \notag
\end{align}
After possibly replacing $S$ by 
its \'{e}tale covering,
 we may assume  that {\it $S$ is affine and  the vector bundle $\Omega_{S/K}$ is free}.
Under this assumption, 
$\Psi_S^*(\omega^{\mr{PGL}}_{S}) - \omega^{\mr{can}}_{S} =0$ if and only if
$\sigma_{S, A}^{T^\vee *}(\Psi_S^*(\omega^{\mr{PGL}}_{S}) - \omega^{\mr{can}}_{S}) =0$ for all
$A \in \Gamma (S, \Omega_{S/K})$.
Thus, in order to complete the proof of Theorem \ref{t05},  it suffices (by (\ref{E70007})) to prove  the equality
\begin{align} \label{e72}
\sigma_S^*(\omega^\mr{PGL}_S)  = 0_S^*(\omega^{\mr{can}}_{S}).
\end{align}

\vspace{3mm}
\subsection{} \label{s37}
On the one hand,  let us consider the right-hand side  of the required equality  (\ref{e72}), i.e., $ 0_S^*(\omega^{\mr{can}}_{S})$.
The differential of the zero section $0_S : S \migi T^\vee_S$ specifies 
a split injection $\mcT_{S/K} \migiincl 0^*_S(\mcT_{T^\vee_S/K})$ of the natural short exact sequence
\begin{align} 
   0 \migi  0_S^*(\mcT_{T^\vee_S/S}) \migi 0^*_S(\mcT_{T^\vee_S/K}) \migi \mcT_{S/K} \migi 0.\notag
   \end{align}
This split injection gives  a decomposition
\begin{align} \label{e74} 
0^*_S(\mcT_{T^\vee_S/K}) \isom   
 \mcT_{S/K} \oplus 0_S^*(\mcT_{T^\vee_S/S})
\isom 
  \mcT_{S/K} \oplus  \Omega_{S/K}
 \isom 
  \mbR^1f_{S*}(\mcT_{C_S/S}) \oplus f_{S*}(\Omega_{C_S/S}^{\otimes 2}),
\end{align}
where the last isomorphism follows from (\ref{eE19}).
The $\mcO_S$-bilinear map on $0^*_S(\mcT_{T^\vee_S/K})$ determined by $\omega_S^\text{can}$ is 
given, via this decomposition, 
by the pairing
$\langle-, -\rangle : \mbR^1f_{S*}(\mcT_{C_S/S}) \times f_{S*}(\Omega_{C_S/S}^{\otimes 2}) \migi \mcO_S$ arising from $\int_{C_S, \Omega^{\otimes 2}_{C_S/S}}$ (cf. (\ref{e16})).
More precisely, this bilinear map may be expressed, via (\ref{e74}), as the map given by assigning
\begin{align} \label{e75}
     ((a, b), (a',b')) \mapsto \langle a, b'\rangle - \langle a', b\rangle  \end{align}
 for local sections $a$, $a' \in  \mbR^1f_{S*}(\mcT_{C_S/S})$ and $b$, $b' \in f_{S*}(\Omega_{C_S/S}^{\otimes 2})$.
 
\vspace{3mm}
\subsection{} \label{s38}
On the other hand, let us  consider the left-hand side  of  (\ref{e72}), i.e., $\sigma_S^*(\omega^\mr{PGL}_S)$.
The differential of the section $\sigma_S \ (=\Psi_S \circ 0_S)  : S \migi \mfS_S$ specifies 
 a split injection $ \mcT_{S/K} \migiincl \sigma_S^*(\mcT_{\mfS_S/K})$  of  the short  exact sequence
\begin{align} 
0 \migi \sigma_S^*(\mcT_{\mfS_S/S}) \migi \sigma_S^*(\mcT_{\mfS_S/K}) \migi \mcT_{S/K} \migi 0.   \notag
\end{align}
If  ${C_S}_{/S}^{^\text{Zzz...}} = (C_S/S, \mcP^\ind =(\mcP_B, \nabla_{\mcP_G}))$ denotes  the ordinary dormant curve classified by $\breve{v}$ (hence  $\gamma^\heartsuit_{\mcP^\ind}$ is an isomorphism), then it follows from  Proposition \ref{p02} that
this split injection corresponds to
a  split injection
\begin{align} \label{eE20}
\mbR^1f_{S*}(\mcT_{C_S/S}) \migiincl  \mbR^1f_{S*}(\mcK^\bullet [\nabla_{\mcP_G}^\mr{ad}])
\end{align}
 of   the short   exact sequence
\begin{align} 
 0 \migi f_{S*}(\Omega_{C_S/S}^{\otimes 2}) \xrightarrow{\gamma_{\mcP^\ind}^\sharp}\mbR^1f_{S*}(\mcK^\bullet [\nabla_{\mcP_G}^\mr{ad} ]) \xrightarrow{\gamma^\flat_{\mcP^\ind}} \mbR^1f_{S*}(\mcT_{C_S/S})  \migi 0. \notag
 \end{align}
Moreover,
if we identify $\mbR^1f_{S*}(\mcT_{C_S/S}) $ with  $\mbR^1f_{S*}(\mr{Ker}(\nabla_{\mcP_G}^\mr{ad}))$ via 
the isomorphism $\gamma^\heartsuit_{\mcP^\ind}$, then 
 it follows from Proposition \ref{p07} that
the injection (\ref{eE20})  coincides with 
 $\gamma_{\mcP^\ind}^\natural : \mbR^1f_{S*}(\mr{Ker}(\nabla_{\mcP_G}^\mr{ad})) \migi \mbR^1f_{S*}(\mcK^\bullet [\nabla_{\mcP_G}^\mr{ad}])$.
Consider the  resulting decomposition 
 \begin{align} \label{e78}
   \mbR^1f_{S*}(\mcK^\bullet [\nabla_{\mcP_G}^\mr{ad}]) \isom  \mbR^1f_{S*}(\mcT_{C_S/S})  \oplus f_{S*}(\Omega_{C_S/S}^{\otimes 2}). 
   \end{align}
Because of the discussion in the previous subsection, the proof of Theorem \ref{t05} is reduced  to the following lemma.

\vspace{3mm}
\ble \label{p40}
 The $\mcO_S$-bilinear map on  $\mbR^1f_{S*}(\mcK^\bullet [\nabla_{\mcP_G}^\mr{ad} ]) $ corresponding to  $\omega_S^\mr{PGL}$  is given, via the decomposition  (\ref{e78}), 
  by the pairing $ \mbR^1f_{S*}(\mcT_{C_S/S})  \times  f_{S*}(\Omega_{C_S/S}^{\otimes 2}) \migi \mcO_S$ arising from $\int_{C_S, \Omega^{\otimes 2}_{C_S/S}}$
  (in the sense of (\ref{e75})).
   \ele
\begin{proof}
The subsheaf $\mr{ad}(\mcP_G)^1 \subseteq \mr{ad}(\mcP_G)$ is isotropic with respect to 
$\kappa_{\mcP^\ind} : \mr{ad}(\mcP_G) \otimes_{\mcO_{C_S}} \mr{ad}(\mcP_G) \migi \mcO_{C_S}$ (cf. (\ref{e52})).
Hence,   $\kappa_{\mcP^\ind}$ induces an $\mcO_{C_S}$-bilinear map
$\overline{\kappa}_{\mcP^\ind} : (\mr{ad}(\mcP_G)/\mr{ad}(\mcP_G)^1) \times \mr{ad}(\mcP_G)^2 \migi \mcO_{C_S}$.
One verifies  from a straightforward calculation that
the diagram
\begin{align}
\xymatrix{ \mcT_{C_S/S} \times \Omega_{C_S/S} \ar[rr]^{\hspace{-10mm} (\overline{\nabla}^\flat)^{-1} \times \overline{\nabla}^\sharp} \ar[rd]_{\langle -, - \rangle} & & (\mr{ad}(\mcP_G)/\mr{ad}(\mcP_G)^1) \times \mr{ad}(\mcP_G)^2  \ar[ld]^{\overline{\kappa}_{\mcP^\ind}}\\
& \mcO_{C_S} &
}\notag
\end{align}
is commutative, where $\langle -, - \rangle$ denotes the natural paring $\mcT_{C_S/S}\times \Omega_{C_S/S} \migi \mcO_{C_S}$.
 Therefore,  the assertion follows from 
 this observation
 together with,
 e.g., the  explicit description of $\mbR^1f_{S*}(\mcK^\bullet [\nabla_{\mcP_G}^\mr{ad}])$ in terms of the \v{C}ech double complex (cf. \S\,\ref{Ws3005}).
\end{proof}
\vspace{3mm}

Consequently, we have proved   Theorem \ref{t05} (i.e., Theorem A), as desired.

\vspace{10mm}
\section{Appendix: Application of Theorem A} \label{s39}\vspace{3mm}

As an application of Theorem \ref{t05}, we construct  certain additional structures on  
${^\circledcirc \mfS^{^\mr{Zzz...}}_{g, K}}$.
Let  $K$ be a field of characteristic $p >2$, 
$X$ a smooth Deligne-Mumford stack  over $K$, and $\omega$  a symplectic structure  on  $X$.
The $2$-form $\omega$ corresponds to a nondegenerate pairing $\mcT_{X/K} \otimes_{\mcO_X} \mcT_{X/K} \migi \mcO_X$ on $\mcT_{X/K}$ and gives an identification $\mcT_{X/K} \isom \mcT_{X/K}^\vee = \Omega_{X/K}$.
Under  this identification,  $\omega$ may be thought of as  a nondegenerate pairing
$\omega^{-1} : \Omega_{X/K} \otimes_{\mcO_X} \Omega_{X/K} \migi \mcO_X$.
Thus, we obtain  a skew-symmetric  $K$-bilinear map
\begin{align} 
\{-, - \}_{\omega} : \mcO_X \times \mcO_X \migi \mcO_X
\notag \end{align}
 defined by 
$\{ f, g \}_{\omega} = \omega^{-1}(df, dg)$.
One verifies from the closedness of  $\omega$ that $\{-, -\}_\omega$ defines  a Poisson bracket.

\vspace{3mm}
\bde \label{d10}
A {\bf restricted structure} on the pair  $(X, \omega)$  is a map $(-)^{[p]} : \mcO_X \migi \mcO_X$ such that
 the triple $(\mcO_X, \{-,-\}_\omega, (-)^{[p]})$ forms a sheaf of restricted Poisson algebras over $K$ (cf. ~\cite{BeKa2}, Definition 1.8,  for the definition of a restricted Poisson algebra).
 \ede
\vspace{3mm}

In addition, we shall recall the definition of a Frobenius-constant quantization (cf.  ~\cite{BeKa1}, Definition 3.3; \cite{BeKa2}, Definition 1.1 and Definition 1.4).
We shall denote by $X^{(1)}$  the Frobenius twist of $X$ over $K$, i.e., the base-change $X \times_{K,  F_K} K$ of the $K$-stack  $X$ via the absolute Frobenius morphism $F_K$ of $K$.
Also, denote by $F: X \migi X^{(1)}$  the relative Frobenius morphism of $X$ over $K$.

\vspace{3mm}
\bde \label{d11} 
\begin{itemize}
\item[(i)]
 Let us consider a pair $(\mcO_X^\hslash, \tau)$ consisting of
\begin{itemize}
\item[$\bullet$] 
 a Zariski sheaf $\mcO_X^\hslash$ of flat $k[[\hslash]]$-algebras on  $X$ complete with respect to the $\hslash$-adic filtration; 
 \item[$\bullet$]
   an isomorphism $\tau : \mcO_X^\hslash /\hslash \cdot  \mcO_X^\hslash \isom \mcO_X$ of sheaves of algebras.
\end{itemize}
Then, we shall say that the pair $(\mcO_X^\hslash, \tau)$  is a {\bf quantization} on the pair $(X, \omega)$ if 
the commutator in $\mcO_X^\hslash$ is equal, via $\tau$, to $\hslash \cdot \{-,-\}_\omega$ mod $\hslash^2 \cdot \mcO_X^\hslash$.
\item[(ii)]
A {\bf Frobenius-constant quantization} on $(X, \omega)$ is a collection of data
\begin{align} 
  \mfO^\hslash_X = (\mcO_X^{\hslash}, \tau, s) \notag
  \end{align} 
consisting of a quantization $(\mcO_X^{\hslash}, \tau)$ of $(X, \omega)$ 
and a morphism $s : \mcO_{X^{(1)}} \migi \mcZ^\hslash \  (\subseteq \mcO_X^\hslash)$, where $\mcZ^\hslash$  denotes the center of $\mcO_X^\hslash$, of sheaves of algebras such that
the  composite 
\begin{align}
\mcO_{X^{(1)}} \xrightarrow{s} \mcZ^\hslash \xrightarrow{\mr{incl.}}\mcO_X^\hslash \xrightarrow{\mr{quot.}} \mcO_X^\hslash/\hslash \cdot  \mcO_X^\hslash \xrightarrow{\tau} \ \mcO_X
\notag \end{align}
 coincides with the morphism $F^* : \mcO_{X^{(1)}} \migiincl \mcO_X$ induced by $F$.
 \end{itemize} 
 \ede
\vspace{3mm}

See ~\cite{BeKa2}, the discussion at the end of \S\,1.2 or Theorem 1.23,  for relationships between  restricted structures on $(X,\omega)$ and  Frobenius-constant quantizations on $(X, \omega)$.

\vspace{3mm}
\begin{exa} \label{E01}
 Let $S$ be a smooth Deligne-Mumford stack over $K$.
 In a natural manner, one may construct  a restricted structure and a Frobenius-constant quantization on the 
cotangent bundle $T^\vee_S$ equipped with the symplectic structure $\omega_S^\text{can}$,  as follows.
\begin{itemize}
\item[(i)]
 If $f$ is 
 a local function   on $T^\vee_S$  lifted from an open substack of $S$, then we set   $f^{[p]} := f^p$.
Also, if $\partial$ is 
 a fiber-wise linear function on $T^\vee_S$ over an open substack of $S$,
 then
  we define  $\partial^{[p]}$  to  be
  the $p$-th iterate of $\partial$. Then, 
  these assignments $(-)\mapsto (-)^{[p]}$ define an endomorphism
  \begin{align} 
(-)^{[p]}: \mcO_{T^\vee_S} \migi \mcO_{T^\vee_S}\notag
\end{align}
 of  $T^\vee_S$.
 By the discussion in (ii) below and  ~\cite{BeKa2}, the discussion at the end of \S\,1.2, the map   $(-)^{[p]}$ forms a restricted structure on $(T^\vee_S, \omega^\mr{can}_S)$.

\item[(ii)]
Assume  that $S$ is  affine. 
 The {\it sheaf of asymptotic differential operators $D^\hslash (S)$} on $S$ (cf. ~\cite{BeKa1}, Example 3.1) is  
the $\hslash$-completion of the $K[\hslash]$-algebra generated by $\Gamma (S, \mcO_S)$ and $\Gamma(S,\mcT_{S/K})$ subject to the  relations
\begin{align}
f_1 *  f_2 = f_1 \cdot f_2, \hspace{4mm}
f_1 * \xi_1 = f_1 \cdot  \xi_1, \hspace{4mm}
\xi_1 * \xi_2 -\xi_2 * \xi_2= \hslash \cdot  [\xi_1, \xi_2], \hspace{4mm}
\xi_1 * f_1 - f_1 *  \xi_1 = \hslash \cdot  \xi_1 (f_1). \notag
\end{align}
for  any local sections  $f_1$,  $f_2 \in \Gamma(S, \mcO_S)$ and $\xi_1$, $\xi_2 \in \Gamma(S, \mcT_{S/K})$, where $*$ denotes the multiplication in $D^\hslash (S)$.
We have a natural isomorphism 
\begin{align}
\tau (S) : D^\hslash (S)/\hslash \cdot D^\hslash (S) \isom  \Gamma(T^\vee_S, \mcO_{T^\vee_S}).\notag
\end{align}
If $Z^\hslash (S)$ denotes the center of $D^\hslash (S)$, then 
 there exists (cf. ~\cite{BeKa1}, the proof of Proposition 3.5) a morphism
\begin{align}
s (S) : \Gamma (T^{\vee (1)}_S, \mcO_{T^{\vee (1)}_S}) \migi Z^\hslash  (S)\notag
\end{align}
  determined, under the identification  $\Gamma (T^{\vee (1)}_S, \mcO_{T^{\vee (1)}_S})=\Gamma (T^\vee_S, \mcO_{T^\vee_S})$ via 
$\mr{id}_{T^\vee_S} \times F^\mr{abs}_K : T^{\vee (1)}_S \isom T^\vee_S$, by 
 \begin{align}
s (S) (f) = f^p,
 \hspace{10mm}
  s (S) (\partial) = \partial^p - \hslash^{p-1} \cdot \partial^{[p]} \notag
 \end{align}
 for any $f \in \Gamma (S, \mcO_S)$,  and $\partial \in \Gamma (S, \mcT_{S/K})$. 
By applying a natural noncommutative localization procedure called {\it Ore localization} (cf. ~\cite{GW}, \S\,9, p.143, Definition),
we obtain, from the triple $(D^\hslash (S), \tau(S), s(S))$,  a Frobenius-constant quantization 
\begin{align} \label{e84}
 (\mcD^\hslash_S, \tau, s)
 \end{align}
 on $(T^\vee_S, \omega_S^\mr{can})$ (cf. ~\cite{BeKa1}, Proposition 3.5).

In general,  let $S$ be 
 an arbitrary   smooth Deligne-Mumford stack.
  Then,  
the Frobenius-constant quantizations constructed above  for various affine schemes having an open immersion into  $S$ may be glued together to  obtain 
a Frobenius-constant quantization
on $(T^\vee_S, \omega_S^\mr{can})$. 
\end{itemize}
\end{exa}
\vspace{3mm}

By applying the discussion in Example \ref{E01} to the case $S = {{^\circledcirc \Dp_{g,K}}}$,
we obtain a restricted structure as well as a Frobenius-constant quantization on $( {{^\circledcirc T^{\vee ^\mr{Zzz...}}_{g,K}}}, \omega^\text{can}_\circledcirc)$.
Such additional structures may be evidently  transported  into $({^\circledcirc \mfS^{^\mr{Zzz...}}_{g, K}}, \omega^\mr{PGL}_\circledcirc)$  via the  isomorphism $\Psi_{g,K}$.
The isomorphism  $\Psi_{g,K}$  preserves   the  symplectic structure as asserted in   Theorem \ref{t05}, so  we have the following  corollary.

\vspace{3mm}
\bco \label{p76}
There exist a  canonical restricted structure and a Frobenius-constant quantization on  
$({^\circledcirc \mfS^{^\mr{Zzz...}}_{g, K}}, \omega^\mr{PGL}_\circledcirc)$.
 \eco

\end{document}